\documentclass[preprints,article,accept,moreauthors,pdftex]{my-mdpi} 

\firstpage{1}
\pubvolume{5}
\issuenum{4}
\articlenumber{261}
\doinum{10.3390/fractalfract5040261}
\pubyear{2021}
\copyrightyear{2021}
\externaleditor{Academic Editor: Vassili Kolokoltsov and Omar Abu Arqub} 
\datereceived{20 October 2021}
\daterevised{19 and 29 November 2021}
\dateaccepted{2 December 2021}
\datepublished{7 December 2021}
\hreflink{https://doi.org/10.3390/\linebreak fractalfract5040261}

\pdfoutput=1

\usepackage{amsmath,mathtools}
\usepackage[labelformat=simple]{subfig}

\usepackage{graphicx}

\newcommand{\Dl}{{}^{\scriptscriptstyle C}_{\scriptscriptstyle 0}\!D^{\alpha}_{\scriptscriptstyle t}}

\Title{Fractional-Order Modelling and Optimal Control of Cholera~Transmission}

\TitleCitation{Fractional-Order Modelling and Optimal Control of Cholera Transmission}

\Author{Silv\'erio Rosa $^{1,}$*$^{,\dagger}$\orcidA{} and Delfim F. M. Torres $^{2,\dagger}$\orcidB{}} 
\AuthorNames{Silv\'erio Rosa and Delfim F. M. Torres}

\AuthorCitation{Rosa, S.; Torres, D.F.M.}

\address{%
$^{1}$ \quad Instituto de Telecomunicações (IT) and Department of Mathematics, 
Universidade da Beira Interior, 6201-001~Covilh\~{a}, Portugal\\
$^{2}$ \quad Center for Research and Development in Mathematics and Applications (CIDMA), 
\mbox{Department of Mathematics}, University of Aveiro, 3810-193 Aveiro, Portugal; delfim@ua.pt}

\corres{Correspondence: rosa@ubi.pt}

\firstnote{These authors contributed equally to this work.}

\abstract{A Caputo-type fractional-order mathematical model for
``meta\-po\-pu\-la\-tion cholera transmission'' was
recently proposed in [Chaos Solitons Fractals 117 (2018), 37--49].
A sensitivity analysis of that model is done here to show
the accuracy relevance of parameter estimation. Then, a fractional
optimal control (FOC) problem is formulated and numerically solved.
A cost-effectiveness analysis is performed to assess the relevance
of studied control measures. Moreover, such analysis allows us to assess the
cost and effectiveness of the control measures during intervention.
We conclude that the FOC system is more effective only in part of the time interval.
For this reason, we propose a system where the derivative order varies along
the time interval, being fractional or classical when more advantageous.
Such variable-order fractional model, that we call a \emph{FractInt} system,
shows to be the most effective in the control of the disease.}

\keyword{cholera;
compartmental mathematical models;
fractional-order optimal control}

\MSC{34A08, 49M05, 92C60}

\begin{document}


\section{Introduction}

Fractional calculus is an old subject that raised as a consequence of a pertinent
question that L'H\^{o}pital asked Leibniz in a letter about
the possible meaning of a derivative of order $1/2$.~Recently,
many researchers have focused their attention in modelling
real-world phenomena using fractional-order derivatives.
The dynamics of those problems have been modelled and studied
by using the concept of fractional-order derivatives. Such problems
appear in, for~example, biology, physics, ecology, engineering,
and various other fields of applied sciences, see, e.g.,~\cite{MR1658022,MR3935222,MR4281822}.

\textls[-5]{Cholera is a gastroenteritis infection, contracted after consuming
an infectious dose or inoculum size of the pathogenic \emph{vibrio cholerae}
\cite{Mukandavire8767}. The~mode of transmission consists of two pathways:
the primary route, where individuals consumes the pathogen
from \emph{vibrio} contaminated water and seafood; the secondary route
being characterised by individuals consuming unhygienic or soiled
food that is infested with pathogenic \emph{vibrios} from an infected person.
This secondary route of transmission is also commonly referred to
as person-to-person contact~\cite{Mukandavire8767}. Cholera infection
has affected many parts of the world. However, its devastating force
has been more pronounced in impoverished communities~\cite{MR3602689,MyID:417}.}

Cholera is  one of the most studied infections in recent years.
Mathematical models are used to study and understand the dynamics
of infection, as well as offer suggestions toward its control, see, e.g.,
\cite{Mukandavire8767,MR4110649,Codeco2001,Capasso1979121,%
njagarah2018spatial,njagarah2015modelling,neilan2010modeling}.
Most models used in the study of Cholera have been based on systems of integer-order
ordinary differential equations. However, those models do not fully account for memory,
as well as non-local properties exhibited by the epidemic system. Non-local behaviour
asserts that the subsequent state of the model depends on both the current and historical states.
Fractional order differential systems have been proposed as more suitable to describe
epidemic dynamics of diseases~\cite{MyID:471,MyID:483,MR4232864}.

Recently, a~fractional order differential system was used in the
study of cholera transmission in adjacent communities~\cite{njagarah2018spatial}.~Individuals 
in adjacent communities often frequent their home range,
an aspect which constitutes memory that is present in such a model.
Here, we start to do a sensitivity  analysis to the fractional-order model
of~\cite{njagarah2018spatial} in order to determine which model parameters
are most influential on the disease dynamics. After~that,
fractional optimal control (FOC) is applied, as~a generalisation
of the classical optimal control system of~\cite{njagarah2015modelling}.
In contrast with other works, for~example~\cite{Rosa}, we could
not find a derivative order for which the FOC system is more effective.
However, we noticed that the FOC system is more effective in part of the time interval.
Hence, we propose a system where the derivative order varies along the interval,
being fractional or classical when more advantageous. Such system,
a variable-order fractional one, is named here a \emph{FractInt} system,
shown to be useful in the control of the~disease.

The paper is organised as follows. In~Section~\ref{2:section}, the~fractional
order model formulation is presented.~Our main results are then given in
Section~\ref{3:mainresults}: sensitivity analysis of the parameters of
the model, taking into account the derivative order (Section~\ref{31:sensitivity});
fractional optimal control of the model (Section~\ref{subsec:control});
numerical simulations and  cost-effectiveness of the fractional model
(Section~\ref{subsec:numres}); and the new variable-order \emph{FractInt} 
system (\mbox{Section~\ref{subsec:fractint}}). We end with Section~\ref{sec:conclusion} 
of~conclusions.


\section{Fractional-Order Cholera~Model}
\label{2:section}

A metapopulation model for cholera transmission is considered,
dividing the population into mutually exclusive distinct
groups and using deterministic continuous transitions between those groups,
also known as states.~The~model describes the dynamics of a population exposed to
infection by the pathogen \emph{vibrio cholerae}.~The~human population is divided
into three compartments: susceptible individuals ($S$);
infectious individuals ($I$); and recovered individuals ($R$).
The Caputo fractional-order system of differential equations
is as follows~\cite{njagarah2018spatial}:
\begin{equation}
\label{Cholera_model:com1}
\left\{
\begin{split}
\Dl S_1 = & \pi_1 +a_2 S_2+\omega R_1-(1-u)\dfrac{\beta_1 B_1 S_1}{K+B_1}\\[1mm]
& -(1-m)\varrho_1 I_1 S_1-(a_1+\mu_1+v)S_1,\\[1mm]
\Dl I_1 = & (1-u)\dfrac{\beta_1 B_1 S_1}{K+B_1}+(1-m)\varrho_1 I_1 S_1+b_2 I_2-Q_1 I_1,\\[1mm]
\Dl R_1 = & vS_1 +\gamma_1 I_1-(\mu_1+\omega+c_1) R_1+c_2 R_2,\\[1mm]
\Dl B_1 = & \sigma_1 I_1-Q_2 B_1,\\[1mm]
\end{split}\right.
\end{equation}
for the first sub-population, and~\begin{equation}
\label{Cholera_model:com2}
\left\{
\begin{split}
\Dl S_2 = & \pi_2+a_1 S_1+\omega R_2-(1-u)\dfrac{\beta_2 B_2 S_2}{K+B_2}\\[1mm]
&-(1-m)\varrho_2 I_2 S_2-(a_2+\mu_2+v) S_2,\\[1mm]
\Dl I_2 = & (1-u)\dfrac{\beta_2 B_2 S_2}{K+B_2}+(1-m)\varrho_2 I_2 S_2+b_1 I_1-Q_3 I_2,\\[1mm]
\Dl R_2 = & vS_2+\gamma_2 I_2-(\mu_2+\omega+c_2) R_2+c_1 R_1,\\[1mm]
\Dl B_2 = & \sigma_2 I_2-Q_4 B_2,
\end{split}\right.
\end{equation}
for the second sub-population, where $\Dl$ denotes the left Caputo
fractional order derivative of order $\alpha$ \cite{MR1658022},
$0 < \alpha \leqslant 1$,
$Q_1=\mu_1+\delta_1+\gamma_1+b_1$,
$Q_2=\mu_p-g_1$, $Q_3=\mu_2+\delta_2+\gamma_2+b_2$,
and $Q_4=\mu_p-g_2$.

We note that the equations of model \eqref{Cholera_model:com1}--\eqref{Cholera_model:com2}
do not have appropriate time dimensions. Indeed, on~the left-hand side the dimension
is (time)$^{-\alpha}$ while on the right-hand side the dimension is (time)$^{-1}$
(see, e.g.,~\cite{MR3808497,MR3928263} for more details).~Therefore, we claim
that the accurate way of writing system \eqref{Cholera_model:com1} is
\begin{equation}
\label{Cholera_model:I}
\left\{\begin{split}
\Dl S_1 = & \pi_1^{\alpha}  +a_2^{\alpha} S_2
+\omega^{\alpha} R_1-(1-u)\dfrac{\beta_1^{\alpha} B_1 S_1}{K+B_1}\\[1mm]
&-(1-m)\varrho_1^{\alpha} I_1 S_1-(a_1^{\alpha}+\mu_1^{\alpha}+v)S_1,\\[1mm]
\Dl I_1 = & (1-u)\dfrac{\beta_1^{\alpha} B_1 S_1}{K+B_1}
+(1-m)\varrho_1^{\alpha} I_1 S_1+b_2^{\alpha} I_2-Q_1 I_1,\\[1mm]
\Dl R_1 = &  vS_1+\gamma_1^{\alpha} I_1-(\mu_1^{\alpha}
+\omega^{\alpha}+c_1^{\alpha}) R_1+c_2^{\alpha} R_2,\\[1mm]
\Dl B_1 = & \sigma_1^{\alpha} I_1-Q_2 B_1,\\[1mm]
\end{split}
\right.
\end{equation}
for the first sub-population, and~\begin{equation}
\label{Cholera_model:II}
\left\{
\begin{split}
\Dl S_2 = & \pi_2^{\alpha}+a_1^{\alpha} S_1
+\omega^{\alpha} R_2-(1-u)\dfrac{\beta_2^{\alpha} B_2 S_2}{K+B_2}\\[1mm]
&-(1-m)\varrho_2^{\alpha} I_2 S_2-(a_2^{\alpha}+\mu_2^{\alpha}+v) S_2,\\[1mm]
\Dl I_2 = & (1-u)\dfrac{\beta_2^{\alpha} B_2 S_2}{K+B_2}
+(1-m)\varrho_2^{\alpha} I_2 S_2+b_1^{\alpha} I_1-Q_3 I_2,\\[1mm]
\Dl R_2 = & vS_2+\gamma_2^{\alpha} I_2-(\mu_2^{\alpha}
+\omega^{\alpha}+c_2^{\alpha}) R_2+c_1^{\alpha} R_1,\\[1mm]
\Dl B_2 = & \:\sigma_2^{\alpha} I_2-Q_4 B_2,
\end{split}\right.
\end{equation}
the accurate way of writing system \eqref{Cholera_model:com2},
with $Q_1=\mu_1^{\alpha}+\delta_1^{\alpha}+\gamma_1^{\alpha}+b_1^{\alpha}$,
$Q_2=\mu_p^{\alpha}-g_1^{\alpha}$, $Q_3=\mu_2^{\alpha}
+\delta_2^{\alpha}+\gamma_2^{\alpha}+b_2^{\alpha}$,
and $Q_4=\mu_p^{\alpha}-g_2^{\alpha}$.


\section{Main~Results}
\label{3:mainresults}

We begin by doing a sensitivity analysis to the parameters of the model,
in order to identify those for which a small perturbation
leads to relevant quantitative changes in the~dynamics.


\subsection{Sensitivity~Analysis}
\label{31:sensitivity}

Two distinct ways to compute the basic reproduction numbers,
${R_0}_1$ and ${R_0}_2$, of~the two sub-populations of the model,
are available in~\cite{njagarah2018spatial,njagarah2015modelling}. To~know
which one is proper, we determine them by using the next-generation
matrix method~\cite{MR1950747}.~We obtain the community specific
reproduction numbers as
\begin{equation}
\label{r0:Ch01_model}
{R_0}_1=\frac{(\pi_1^{\alpha}(\mu_2^{\alpha}+a_2^{\alpha}+v)
+a_2^{\alpha}\pi_2^{\alpha})((1-u)\beta_1^{\alpha}\sigma_1^{\alpha}
+(1-m)Q_2\varrho_1^{\alpha} K)}{Q_1 Q_2(\mu_1^{\alpha}
+a_1^{\alpha}+v)(\mu_2^{\alpha}+a_2^{\alpha}+v)(1-\Phi_1)K},
\end{equation}
for the first community, and~\begin{equation}
\label{r0:Ch02_model}
{R_0}_2=\frac{(\pi_2^{\alpha}(\mu_1^{\alpha}+a_1^{\alpha}+v)
+a_1^{\alpha}\pi_1^{\alpha})((1-u)\beta_2^{\alpha}\sigma_2^{\alpha}
+(1-m)Q_4\varrho_2^{\alpha} K)}{Q_3 Q_4(\mu_1^{\alpha}
+a_1^{\alpha}+v)(\mu_2^{\alpha}+a_2^{\alpha}+v)(1-\Phi_1)K},
\end{equation}
for the second community, where
\[
\Phi_1=\frac{a_1^{\alpha} a_2^{\alpha}}{(\mu_1^{\alpha}
+a_1^{\alpha}+v)(\mu_2^{\alpha}+a_2^{\alpha}+v)}.
\]

\begin{figure}[H]
\subfloat{%
\resizebox*{7cm}{!}{\includegraphics{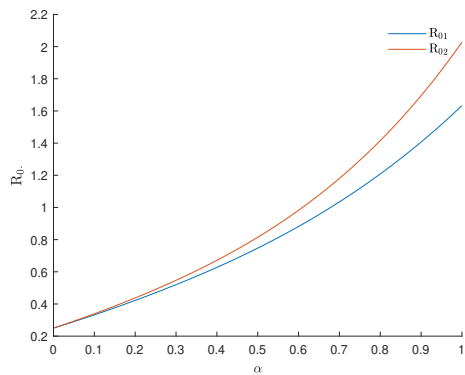}}}
\caption{Impact of the variation of derivative order, $\alpha$,
in the basic reproduction numbers of both communities \eqref{Cholera_model:I}
and \eqref{Cholera_model:II}, respectively, ${R_0}_1$ and ${R_0}_2$.}
\label{fig:R0_varalpha}
\end{figure}

The values of used parameters are presented in Table~\ref{tab:param},
which were taken from~\cite{njagarah2018spatial}, with~exception of $\pi_1$,
$\pi_2$ and $\varrho_2$. The~first two parameters are equal and are defined as
$$
\pi_1=\pi_2= 1.08\times10^{-4},
$$
which is bigger than the values proposed in~\cite{njagarah2018spatial},
in order to  ensure an endemic scenario (${R_0}_1>1$ and ${R_0}_2>1$).
Indeed, it is the presence of an endemic situation that motivates us,
in Section~\ref{subsec:control}, to~apply optimal control theory
to tackle the cholera problem. The~latter parameter, $\varrho_2$,
is also changed with respect to~\cite{njagarah2018spatial} and is defined
as $\varrho_2=0.1875$. This particular value makes
specific reproduction numbers of the two communities,
${R_0}_1$ and ${R_0}_2$, clearly different from each other
as we can see in Figure~\ref{fig:R0_varalpha}. This means that
the two communities are distinct and that they could
not be considered as one unique community.
\begin{specialtable}[H]
\caption{Values of model's~parameters.}
\begin{tabular*}{\hsize}{@{}@{\extracolsep{\fill}}ccc@{}}
\toprule
\textbf{Parameter}  & \textbf{Value} & \textbf{Source}\\
\midrule
$\beta_1$  & 0.00125&\cite{njagarah2018spatial}\\
$\beta_2$ & 0.0125&\cite{njagarah2018spatial}\\
$K$ & $10^6$ & \cite{Codeco2001}\\
$\mu_1,\mu_2$  & $8.4\!\times\! 10^{-5}$ & \cite{StatsSA,jamison2006disease}\\
$\delta_1$ &  0.0125 & \cite{neilan2010modeling,KUMATE1998217}\\
$\delta_2$ &  0.045 & \cite{neilan2010modeling,KUMATE1998217}\\
$\gamma_1$ &  0.045 & \cite{neilan2010modeling,king2008inapparent}\\
$\gamma_2$ &  0.035 & \cite{neilan2010modeling,king2008inapparent}\\
$\mu_p$ &  1.06 & \cite{Codeco2001,hartley2005hyperinfectivity,mukandavire2011modelling,munro1996fate}\\
$g_1,g_2$ &  0.73 & \cite{njagarah2018spatial}\\
$\varrho_1$ &  0.102 & -- \\
$\sigma_1,\sigma_2$  & 50 & \cite{njagarah2018spatial}\\
\bottomrule
\end{tabular*}
\label{tab:param}
\end{specialtable}

We start by considering that rates $u$, $v$, and $m$
are all zero in the sensitivity analysis,
unless specified~otherwise.

\begin{Definition}[See~\cite{chitnis2008determining,rodrigues2013sensitivity}]
\label{def:sentInd}
The normalised forward sensitivity index of ${R_0}_i$, which is differentiable
with respect to a given parameter $p$, is defined by
\begin{equation}
\Upsilon_p^{{R_0}_i}=\frac{\partial {R_0}_i}{\partial p}
\frac{p}{{R_0}_i},\quad {\scriptstyle i=1,2}.\label{eq:sentInd}
\end{equation}
\end{Definition}

Table~\ref{tab:sensitivity1} presents the values of
the sensitivity index of the parameters of the model, 
obtained by the normalised sensitivity index \eqref{eq:sentInd}, 
for the classical case ($\alpha=1$) of connected communities.~These 
values have a meaning.~For~instance, $\Upsilon_{\varrho_1}^{{R_0}_1}=+0.999$
means that increasing (decreasing) $\varrho_1$ by a given percentage increases
(decreases) always ${R_0}_1$ by nearly that same percentage.
Sensitive parameters should be carefully evaluated,
once a small perturbation in such parameter leads to significant
quantitative changes. On~the other hand, the~estimation of
a parameter with a small value for the sensitivity index
does not require as much attention to evaluate,
because a small perturbation in that parameter leads
to small adjustments~\cite{mikucki2012sensitivity}.
According with Table~\ref{tab:sensitivity1},
we should pay special attention to the estimation of sensitive
parameters $\varrho_1$ and $\varrho_2$. In~contrast, the~estimation 
of $K$, $\mu_p$, $\beta_1$, $\beta_2$, $\sigma_1$, $\sigma_2$,
$g_1$, and $g_2$ do not require as much attention because of its low sensitivity.
The missing parameters are those whose index value is~zero.

\begin{specialtable}[H]
\caption{Sensitivity of ${R_0}_1$ (\emph{top}) and ${R_0}_2$ (\emph{bottom}),
evaluated for the parameter values given in Table~\ref{tab:param} with $\alpha=1$.}
\label{tab:sensitivity1}
\begin{tabular*}{\hsize}{@{}@{\extracolsep{\fill}}cccccc@{}}
\toprule
Parameter & $\Upsilon_{\cdot}^{{R_0}_1}$ &Parameter & $\Upsilon_{\cdot}^{{R_0}_1}$
&Parameter & $\Upsilon_{\cdot}^{{R_0}_1}$\\[1mm] \midrule
$\pi_1,\pi_2$ &0.500 & $\mu_2,a_1$ & $-0.454$& $\beta_1,\sigma_1$ &$2\times10^{-6}$\\
$\varrho_1$ & 0.999 & $K$ & $-2\times10^{-6}$ &$\mu_1$ &$-0.547$\\
$a_2$ & 0.454 & $b_ 1$ & $-0.343$ & $\delta_1$ & $-0.143$\\
$\gamma_1$ &$-0.514$ & $\mu_p$ & $-6\times10^{-6}$ &$g_1$ &$4\times10^{-6}$\\
\midrule
Parameter & $\Upsilon_{\cdot}^{{R_0}_2}$ & Parameter & $\Upsilon_{\cdot}^{{R_0}_2}$
& Parameter & $\Upsilon_{\cdot}^{{R_0}_2}$\\[1mm] \midrule
$\pi_1,\pi_2$ & 0.500 & $\mu_2$ & $-0.456$ & $\beta_2,\sigma_2$ & $1\times 10^{-5}$\\
$\varrho_2$ & 0.999 & $K$  & $-1\times 10^{-5}$ & $\mu_1,a_2$ & $-0.545$\\
$a_1$ & 0.545 & $b_2$ & $-0.259$ & $\delta_2$ & $-0.416$\\
$\gamma_2$ & $-0.324$ & $\mu_p$ & $-3\times 10^{-5}$ & $g_2$ & $2\times 10^{-5}$\\
\bottomrule
\end{tabular*}
\end{specialtable}

\begin{figure}[H]
\subfloat[]{%
\resizebox*{6.1cm}{!}{\includegraphics{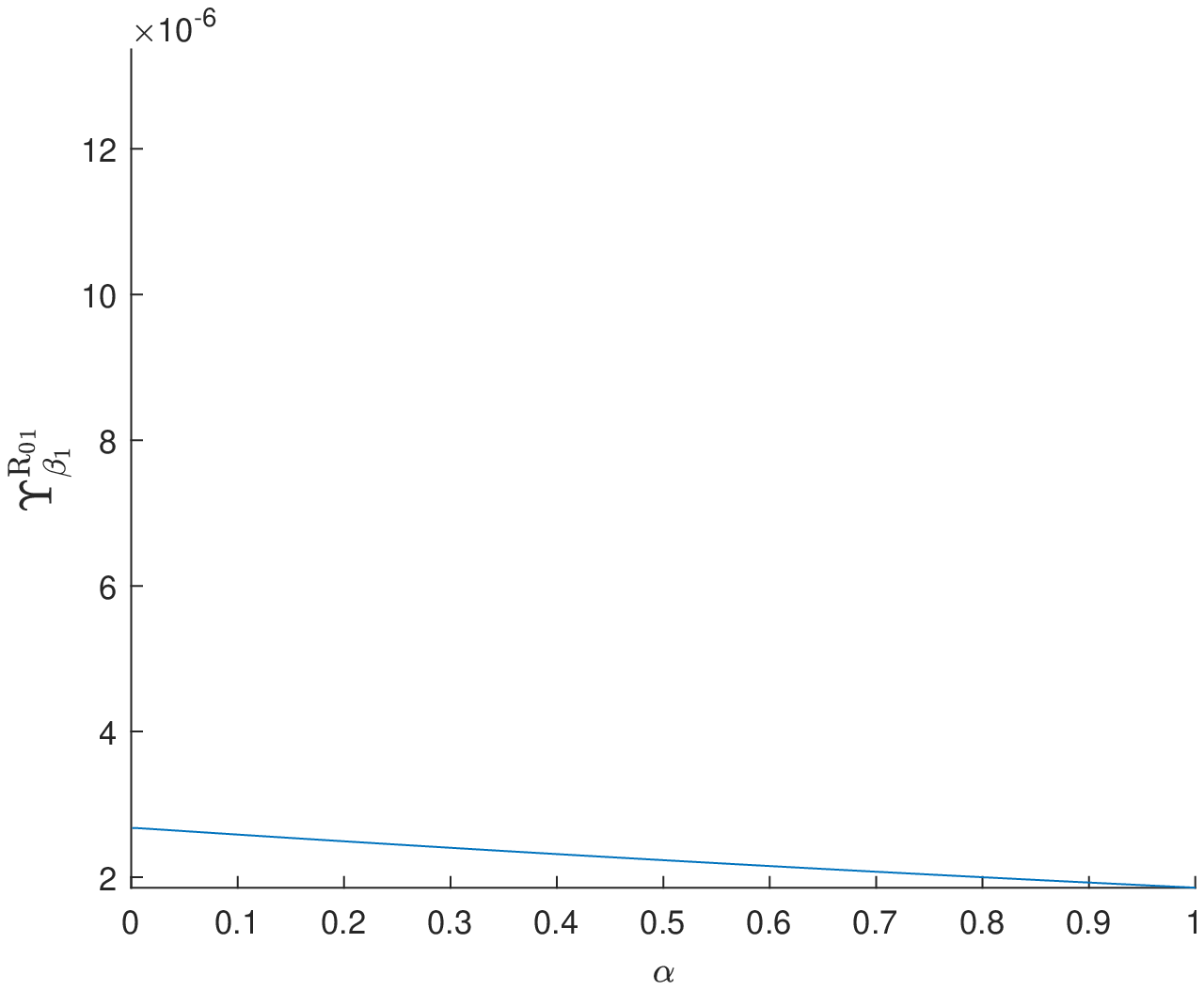}}
\label{fig:sensitivity_beta1}}\hspace{3pt}
\subfloat[]{%
\resizebox*{6.1cm}{!}{\includegraphics{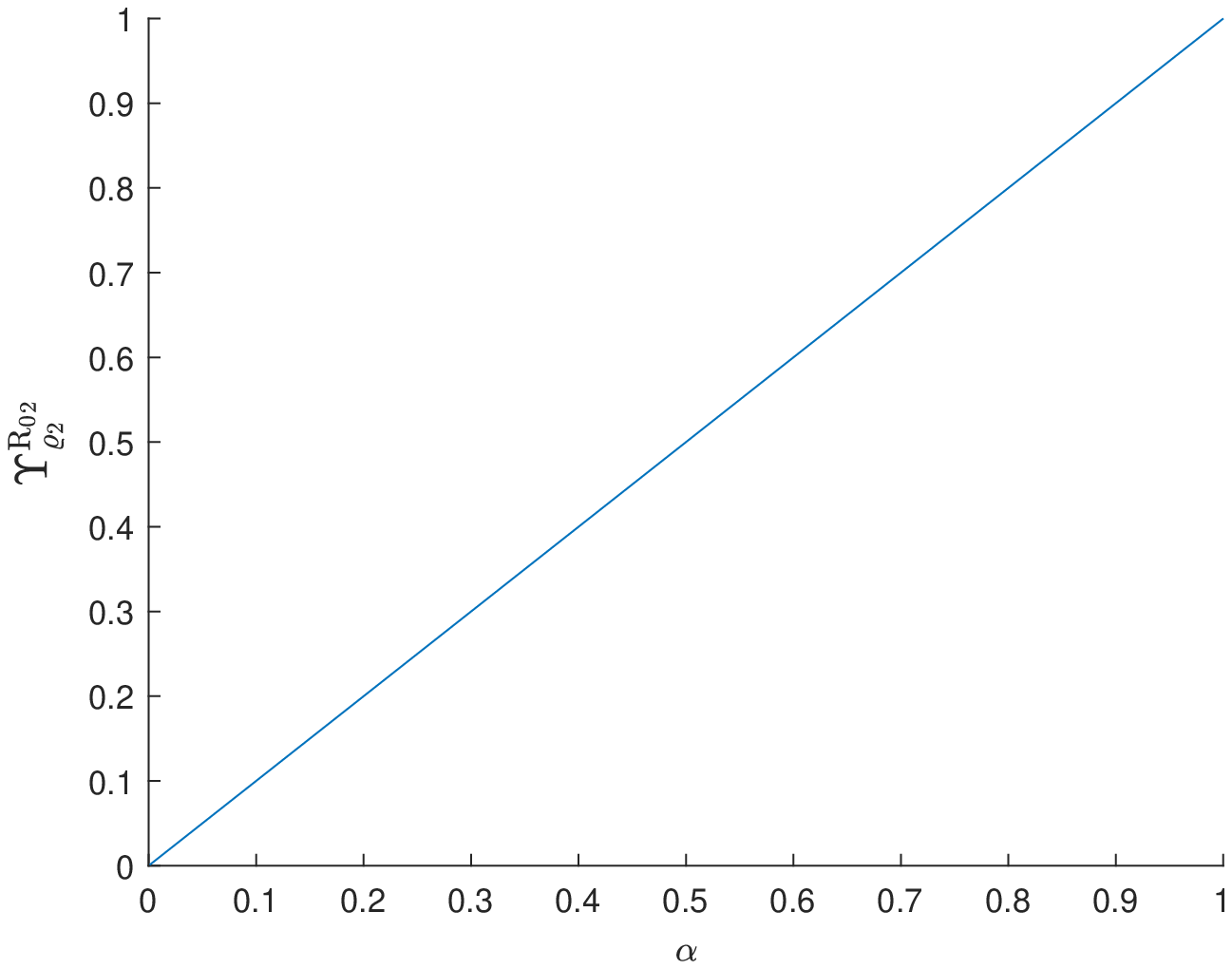}}
\label{fig:sensitivity_varrho2}}
\caption{(\textbf{a})~Evolution of the sensitivity index of parameter $\beta_1$,
evaluated for $\alpha$, with~respect to ${R_0}_1$ (first community); 
(\textbf{b})~Evolution of the sensitivity index of parameter $\varrho_2$, evaluated for $\alpha$,
with respect to ${R_0}_2$ (second community).}
\label{fig:sensitivity}
\end{figure}

For the fractional model, the~sensitivity index depends on the derivative order $\alpha$.
We can see this in Figure~\ref{fig:sensitivity}, where: (a) the impact of variation
of $\alpha$ in the sensitivity index of $\beta_1$ is displayed for the first community;
(b) the impact of variation of $\alpha$ in the sensitivity index of $\varrho_2$
is exhibited for the second community. The~graphics of the other parameters are not shown
because they exhibited similar behaviours. The~parameters whose index value
for $\alpha=1$ is close to zero, as~$\beta_1$, do not vary much if we consider lower values
of $\alpha$, as~we see in Figure~\ref{fig:sensitivity}a. On~the other hand,
parameters whose index value  in Table~\ref{tab:sensitivity1} are not as close to zero
as the previous one, as~$\varrho_2$, vary significantly if we consider lower values of
$\alpha$, as~we see in Figure~\ref{fig:sensitivity}b,
and their sensitivity decreases with the decrease in $\alpha$.

\begin{figure}[H]
\subfloat{%
\resizebox*{7cm}{!}{\includegraphics{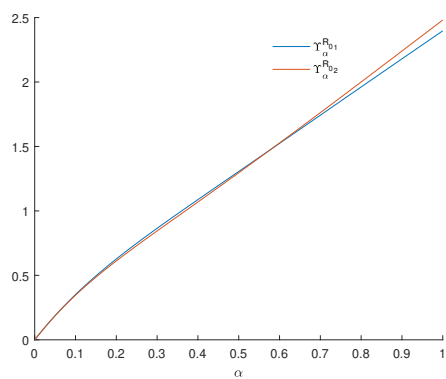}}}
\caption{Evolution of the sensitivity index for the basic reproduction numbers
of both communities, $R{_0}_1$ and ${R_0}_2$, with~the variation
of the derivative order, $\alpha$.}
\label{fig:sensitivity_BRN}
\end{figure}

Figure~\ref{fig:sensitivity_BRN} presents the evolution
of the sensitivity index for the specific reproduction numbers
of the two communities, ${R_0}_1$ and ${R_0}_2$, with~the variation
of the derivative order. We see that the evolution of the sensitivity index
is analogous for the two reproduction~numbers.

\subsection{Fractional Optimal Control of the~Model}
\label{subsec:control}

Modelling dynamic control systems optimally is a very important issue 
in applied sciences and engineering~\cite{MR4271051}.
In this section, our aim is to minimise the number of cholera infected
individuals and, simultaneously, to~reduce the  associated cost. This is
achieved trough: (i) the use of vaccination  into communities, as~an effective
time-dependent measure-control $v(t)$; (ii) the use of clean treated water, 
a~preventive measure control $u(t)$; (iii)  and the implementation of proper
hygiene, another preventive measure control $m(t)$, in~order to control
person-to-person contact. Thus, we consider the following fractional optimal
control problem:
\begin{align}
\label{cost-functional}
\begin{split}
\min ~\mathcal{J}(I_1, I_2,u,v,&m)
=\int_0^{t_f} \left(k_1I_1+ k_2 I_2+k_3 u^2+k_4 v^2+k_5 m^2\right) ~dt\\
\end{split}
\end{align}
subject to system \eqref{Cholera_model:I}--\eqref{Cholera_model:II}
with given initial conditions
\begin{equation}
\label{ocp:ic}
S_1(0)=S_{1_0}\geqslant 0,I_1(0)=I_{1_0}\geqslant 0,\ldots ,
R_2(0)=R_{2_0}\geqslant 0,B_2(0)=B_{2_0}\geqslant 0.
\end{equation}

Note that we are using a quadratic cost functional on the controls,
as an approximation to the real non-linear functional. Indeed,
as in~\cite{MyID:482}, our optimal control problem depends
on the assumption that the cost takes a non-linear form.
The parameters $0<k_1,k_2,k_3,k_4,k_5 <+\infty$ are positive weights and
$t_f$ is the duration of the control program. In addition, $k_3, k_4$, and $k_5$
represent the costs of applying controls efforts $u$, $v$, and $m$, respectively.
The set of admissible control functions is
\begin{equation}
\label{omega:set}
\mathcal{U}=\left\lbrace (u(\cdot),v(\cdot),m(\cdot))
\in L^{\infty}(0,t_{f}) : 0 \leqslant u
\leqslant {u}_{\max},0 \leqslant v
\leqslant {v}_{\max},0 \leqslant m
\leqslant {m}_{\max}
\right\rbrace.
\end{equation}

The Pontryagin maximum principle (PMP) for fractional optimal control
is used to solve the problem~\cite{MR3443073,MR4116679}.~The Hamiltonian 
of the resulting optimal control problem is defined as
\vspace{12pt}
\begin{equation}
\begin{split}
H&= k_1I_1+ k_2 I_2+k_3 u^2+k_4 v^2+k_5 m^2+\xi_1\bigg(\pi_1^{\alpha}  
+a_2^{\alpha} S_2+\omega^{\alpha} R_1\\
&-(1-u)\dfrac{\beta_1^{\alpha} B_1 S_1}{K+B_1}-(1-m)\varrho_1^{\alpha} I_1
S_1-(a_1^{\alpha}+\mu_1^{\alpha}+v)S_1\bigg)\\[1mm]
& +\xi_2\left((1-u)\dfrac{\beta_1^{\alpha} B_1 S_1}{K+B_1}
+(1-m)\varrho_1^{\alpha} I_1 S_1
+b_2^{\alpha} I_2-Q_1 I_1\right)\\[1mm]
&  +\xi_3\left(vS_1+\gamma_1^{\alpha} I_1
-(\mu_1^{\alpha}+\omega^{\alpha}+c_1^{\alpha}) R_1+c_2^{\alpha} R_2\right)
+\xi_4\left(\sigma_1^{\alpha} I_1-Q_2 B_1\right)\\[1mm]
& +\xi_5\bigg(\pi_2^{\alpha}+a_1^{\alpha} S_1
+\omega^{\alpha} R_2-(1-u)\dfrac{\beta_2^{\alpha} B_2 S_2}{K+B_2}\\[1mm]
&-(1-m)\varrho_2^{\alpha} I_2 S_2-(a_2^{\alpha}
+\mu_2^{\alpha}+v) S_2\bigg)\\[1mm]
& +\xi_6\left((1-u)\dfrac{\beta_2^{\alpha} B_2 S_2}{K+B_2}
+(1-m)\varrho_2^{\alpha} I_2 S_2
+b_1^{\alpha} I_1-Q_3 I_2\right)\\[1mm]
& +\xi_7\left(vS_2+\gamma_2^{\alpha} I_2-(\mu_2^{\alpha}
+\omega^{\alpha}+c_2^{\alpha}) R_2+c_1^{\alpha} R_1\right)
+\xi_8\left(\sigma_2^{\alpha} I_2-Q_4 B_2\right)
\end{split}
\end{equation}
and the adjoint system asserts that the co-state variables
$\xi_i(t)$, $i = 1, \ldots, 8$, verify
\begin{align}
\label{eq:co_states_fr2}
\begin{split}
\Dl &\xi_1(t') = \varrho_1^{\alpha} I_1 (m-1) (\xi_1 - \xi_2) +a_1^{\alpha} (\xi_5-\xi_1)\\
&+ \frac{B_1 \beta_1^{\alpha} (\xi_1 - \xi_2) (u-1)-B_1 \mu_1^{\alpha} \xi_1
- K \mu_1^{\alpha} \xi_1 }{B_1 + K} + (\xi_3-\xi_1) v,\\[1.2mm]
\Dl &\xi_2(t')= k_1 + \gamma_1^{\alpha} \xi_3 + b_1^{\alpha} \xi_6
+ \varrho_1^{\alpha} (m-1) \xi_1 S_1 -\xi_2 (Q_1 + \varrho_1^{\alpha} (m-1) S_1)\\
& + \xi_4 \sigma_1^{\alpha},\\[1.2mm]
\Dl &\xi_3(t')  =  \omega^{\alpha} (\xi_1-\xi_3)
- (\mu_1^{\alpha} + c_1^{\alpha}) \xi_3+c_1^{\alpha} \xi_7,\\[1.2mm]
\Dl &\xi_4(t') =\frac{K (-K \xi_4 Q_2 + \beta_1^{\alpha} (\xi_1 -\xi_2) S_1 (u-1))
-B_1^2 \xi_4 Q_2 - 2 B_1 K \xi_4 Q_2}{(B_1 + K)^2},\\
\Dl &\xi_5(t') = a_2^{\alpha} (\xi_1 - \xi_5)
+ \varrho_2^{\alpha} I_2 (m-1) (\xi_5 - \xi_6)\\
& + \frac{B_2 \beta_2^{\alpha} (\xi_5 - \xi_6) (u-1)-B_2 \mu_2^{\alpha} \xi_5
- K \mu_2^{\alpha} \xi_5}{B_2 + K} + (\xi_7-\xi_5) v,\\
\Dl &\xi_6(t') =k_2 + b_2^{\alpha} \xi_2 + \gamma_2^{\alpha} \xi_7
+ \varrho_2^{\alpha} (m-1)\xi_5 S_2
- \xi_6 (Q_3 + \varrho_2^{\alpha} (m-1) S_2)\\
& + \xi_8 \sigma_2^{\alpha},\\
\Dl &\xi_7(t') =\omega^{\alpha} \xi_5 +c_2^{\alpha} \xi_3
- (\mu_2^{\alpha} + \omega^{\alpha}+c_2^{\alpha}) \xi_7,\\
\Dl &\xi_8(t') =\frac{K (-K \xi_8 Q_4 + \beta_2^{\alpha} (\xi_5 - \xi_6)
S_2 (u-1))-B_2^2 \xi_8 Q_4 - 2 B_2 K \xi_8 Q_4}{(B_2 + K)^2},
\end{split}
\end{align}
with $t'=t_f-t$. In~turn, the~optimality conditions of PMP establish
that the optimal controls $u^*$, $v^*$, and $m^*$ are defined by
\begin{equation}
\label{eq:ext:cont}
\begin{split}
u^*&=\min\left\{\max\left\{0,\frac{B_1 \beta_1^{\alpha} (\xi_2-\xi_1) S_1}{2k_3(B_1 + K)}
+ \frac{B_2 \beta_2^{\alpha} (\xi_6-\xi_5) S_2}{2k_3(B_2 + K)} \right\},{u}_{\max}\right\},\\[2mm]
v^*&=\min\left\{\max\left\{0,\frac{S_1(\xi_1  - \xi_3)
+ S_2(\xi_5  - \xi_7)}{2 k_4}  \right\},{v}_{\max}\right\},\\[2mm]
m^*&=\min\left\{\max\left\{0,\frac{\varrho_1^{\alpha} I_1 (\xi_2-\xi_1) S_1
+ \varrho_2^{\alpha} I_2 (\xi_6-\xi_5) S_2}{2 k_5}  \right\},{m}_{\max}\right\}.
\end{split}
\end{equation}
In addition, the~following transversality conditions hold:
\begin{equation}
\label{adjoint:cond_ini}
\xi_i(t_f)=0,
\quad i=1,\ldots,8.
\end{equation}


\subsection{Numerical Results and Cost-Effectiveness~Analysis}
\label{subsec:numres}

\begin{figure}[H]
\subfloat{%
\resizebox*{6.2cm}{!}{\includegraphics{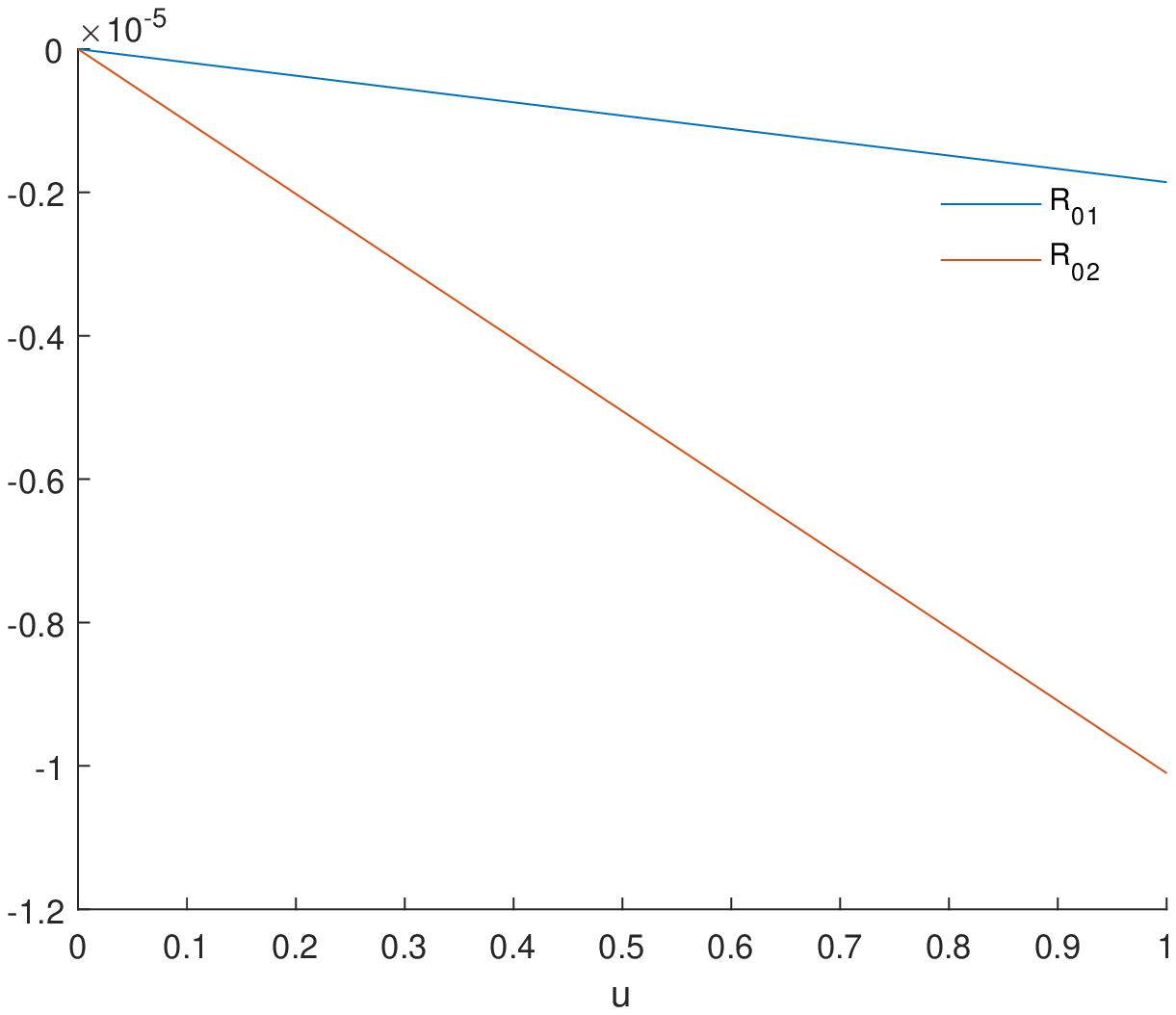}}}\hspace{4pt}
\subfloat{%
\resizebox*{6.2cm}{!}{\includegraphics{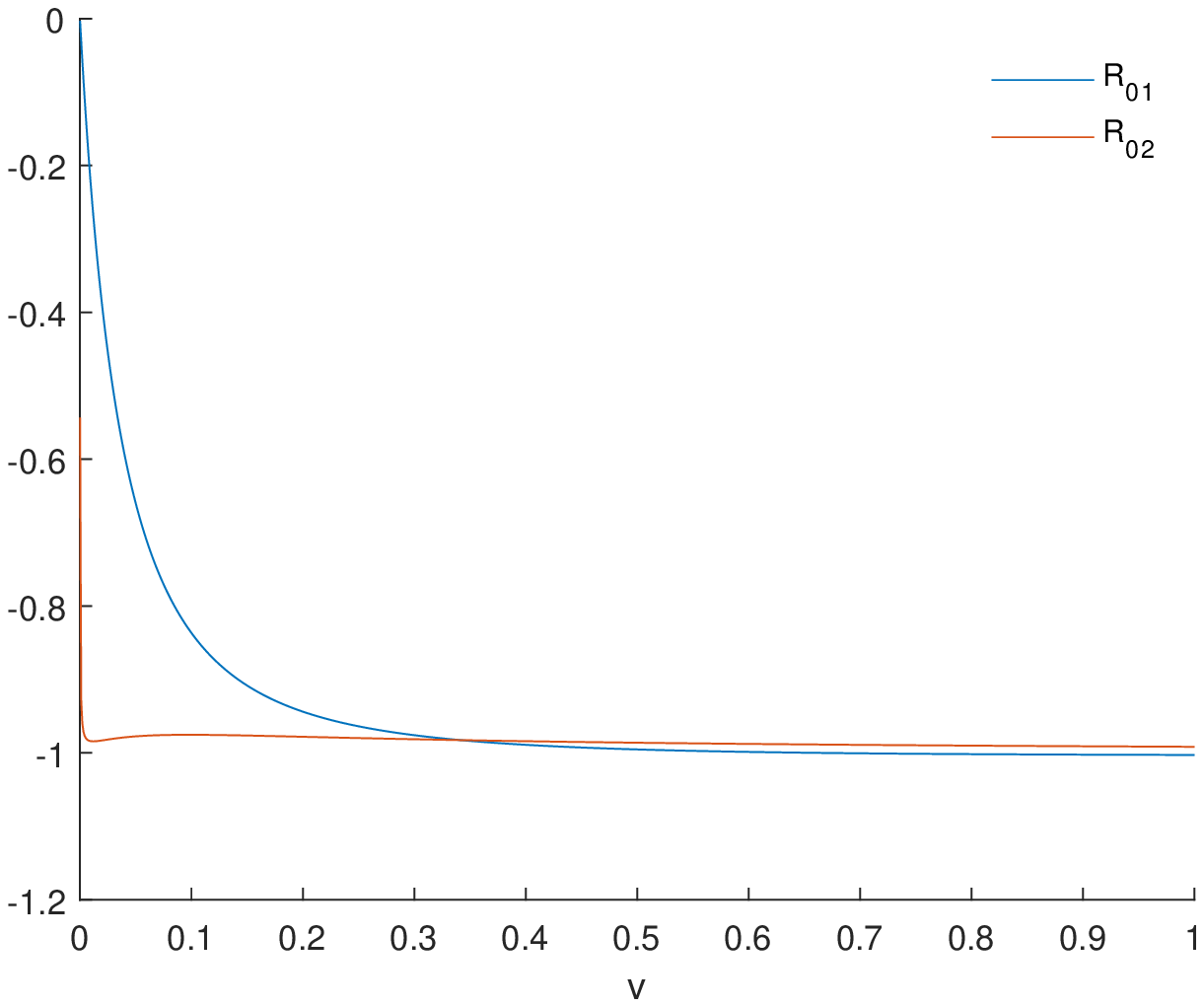}}}\hspace{4pt}
\subfloat{%
\resizebox*{6.2cm}{!}{\includegraphics{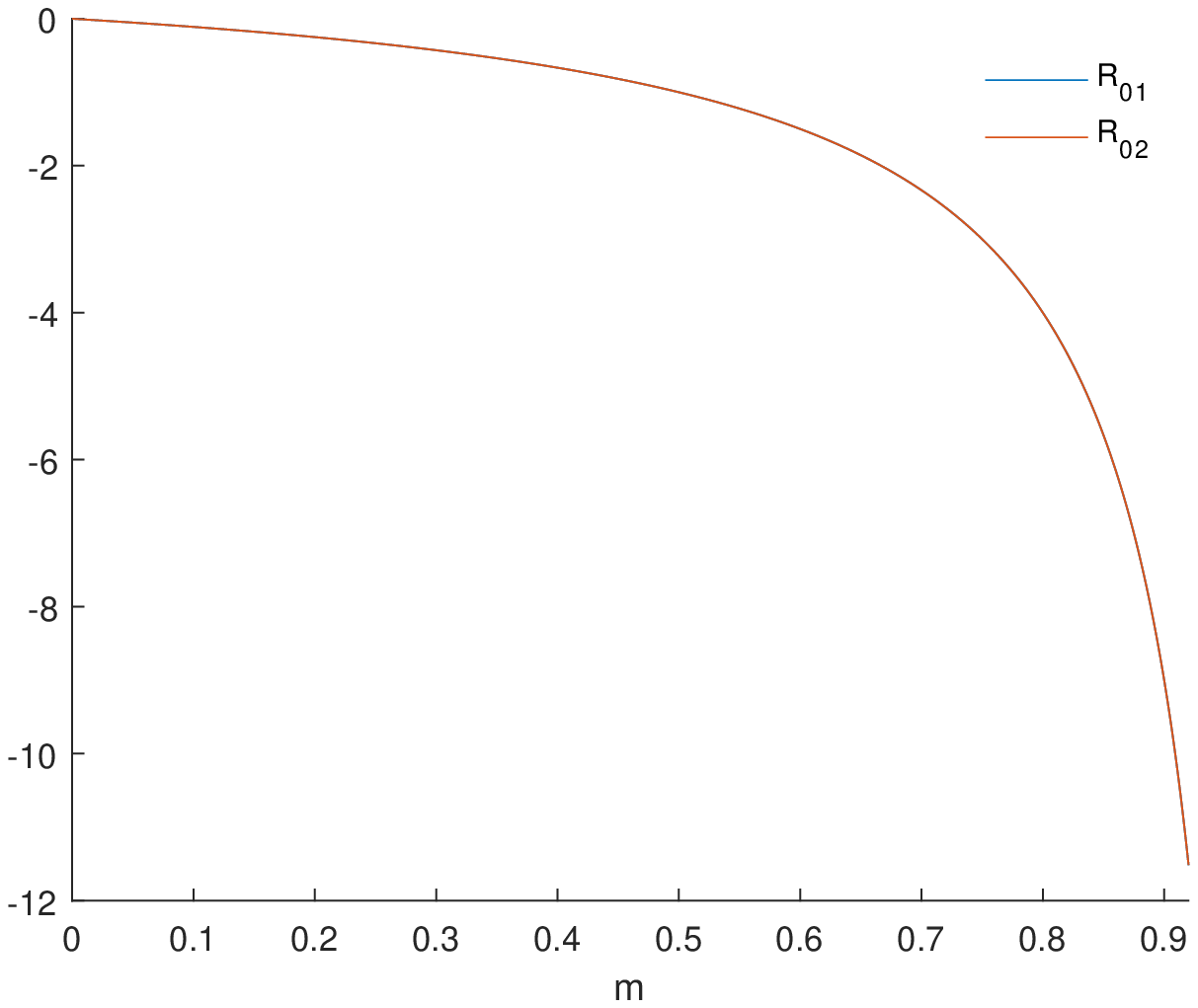}}}
\caption{Sensitivity index of the basic reproduction numbers with respect
to the control variables $u$ (\textbf{top left}), $v$ (\textbf{top right}),
and $m$ (\textbf{bottom}).}
\label{fig:Effic_controls}
\end{figure}

We start by calculating the relevance of the three measures used in the control of
the disease. We do it by using the sensitivity index, presented in Definition~\ref{def:sentInd},
as  proposed in~\cite{2019BpRL...14...27P}. In~this case, the~sensitivity indices
are presented as functions of the control parameters in Figure~\ref{fig:Effic_controls},
using the parametric values from Table~\ref{tab:param} and considering the classical model
(that is, $\alpha=1$). The~resulting graphics show that: (a)~the curve of vaccination
is the one that most rapidly moves away from zero, meaning that the vaccination program
has a big impact for small rates of application; (b)~proper hygiene measures are also
important in the control of cholera, having a bigger impact for greater rates
of application of it, being the only control that has precisely the same impact
in both communities; (c) the domestic water treatment is useless in the control
of cholera transmission, when used simultaneously with the two previous controls.
So, in~what follows, the~third control, variable $u$, is~ignored.

One can deal with many different problems arising in different fields of sciences 
and engineering by applying some appropriate discretisation~\cite{Ed4}.
Here, the Pontryagin Maximum Principle is used to numerically solve the optimal 
control problem, as~discussed in Section~\ref{subsec:control}, both in classical 
and fractional cases, using the predict-evaluate-correct-evaluate (PECE) method of
Adams--Basforth--Moulton~\cite{diethelm2005algorithms}, coded by us in MATLAB.
Firstly, we solve system \eqref{Cholera_model:I}--\eqref{Cholera_model:II}
by the PECE procedure, with~the initial values for the state variables as given
in Table~\ref{tab:solinit} and a guess for the controls over the time interval $[0,t_f]$,
that way obtaining the values of the state variables. Similarly to~\cite{Rosa},
a change of variable is applied to the adjoint system and to the transversality conditions,
obtaining the  fractional initial value problem \eqref{eq:co_states_fr2}--\eqref{adjoint:cond_ini}.
Such IVP is also solved with the PECE algorithm, and~the values of the co-state variables
$\xi_i$, $i=1,\ldots,8$, are determined. The~controls are then updated by a convex combination
of the controls of previous iteration and the current values, computed according
to \eqref{eq:ext:cont}. This procedure is repeated iteratively until the values
of all the variables and the values of the controls are almost identical
to the ones of the previous iteration. The~solutions of the classical model
were successfully confirmed by a classical forward--backward scheme,
also implemented by us in~MATLAB.

In the numerical experiments, we consider weights $k_1=4$,
$k_2=2.4$, $k_3=1.6$, and $k_4=k_5=1$. These values have the same
relations between them than the homonyms weights in~\cite{njagarah2015modelling},
with the numerical advantage of being closer to the value one.
We also use $v_{\max} =m_{\max}= 1$, while the other parameters
are fixed according to Table~\ref{tab:param}.

According with Figure~\ref{fig:R0_varalpha}, ${R_0}_1>1$  and ${R_0}_2>1$
for $\alpha\geqslant 0.68$ (endemic scenario). In~order to be able to compare
the FOCP results for several derivative orders, we consider the  initial conditions,
given by Table~\ref{tab:solinit}, which correspond to the non-trivial endemic equilibrium
for system \eqref{Cholera_model:I}--\eqref{Cholera_model:II} for the classical
cholera model ($\alpha=1$).

\begin{specialtable}[H]
\caption{Initial conditions for the fractional optimal control problem
of Section~\ref{subsec:control} with parameters given by
Table~\ref{tab:param}, corresponding to the endemic equilibrium
of cholera model \eqref{Cholera_model:I}--\eqref{Cholera_model:II}
with classical derivative order.}
\label{tab:solinit}
\begin{tabular*}{\hsize}{@{}@{\extracolsep{\fill}}cccccccc@{}}
\toprule
\boldmath{$S_1(0)$} &\boldmath{ $I_1(0)$} & \boldmath{$R_1(0)$} 
& \boldmath{$B_1(0)$} & \boldmath{$S_2(0)$}&\boldmath{$I_2(0)$}
&\boldmath{$R_2(0)$}& \boldmath{$B_2(0)$} \\[1mm] \midrule
0.53144 & 0.001997 & 0.01028 & 0.30254& 0.44222 & 0.002380 & 0.01082
& 0.36065\\ \bottomrule
\end{tabular*}
\end{specialtable}

\begin{figure}[H]
\subfloat[]{%
\resizebox*{6.1cm}{!}{\includegraphics{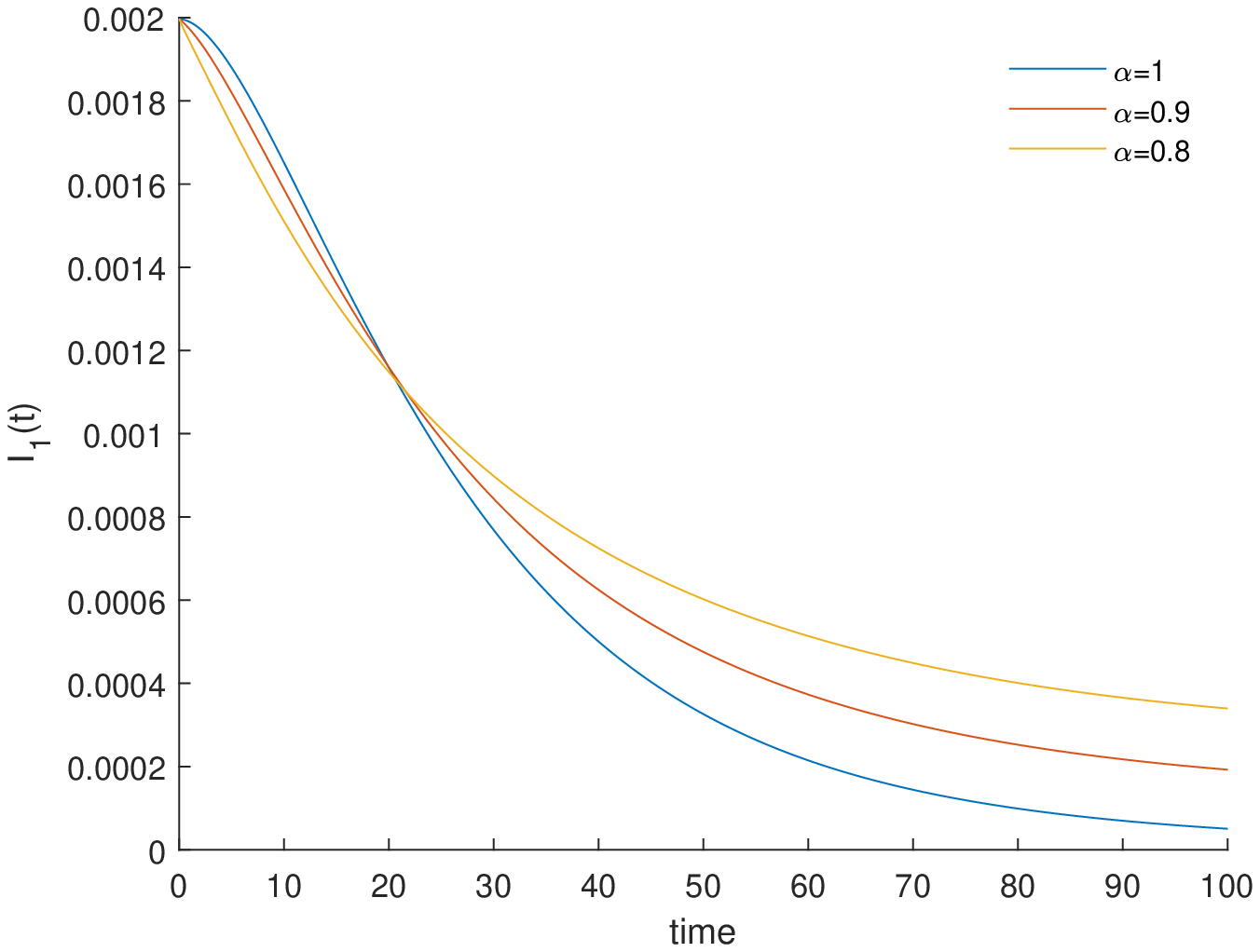}}}\hspace{4pt}
\subfloat[]{%
\resizebox*{6.1cm}{!}{\includegraphics{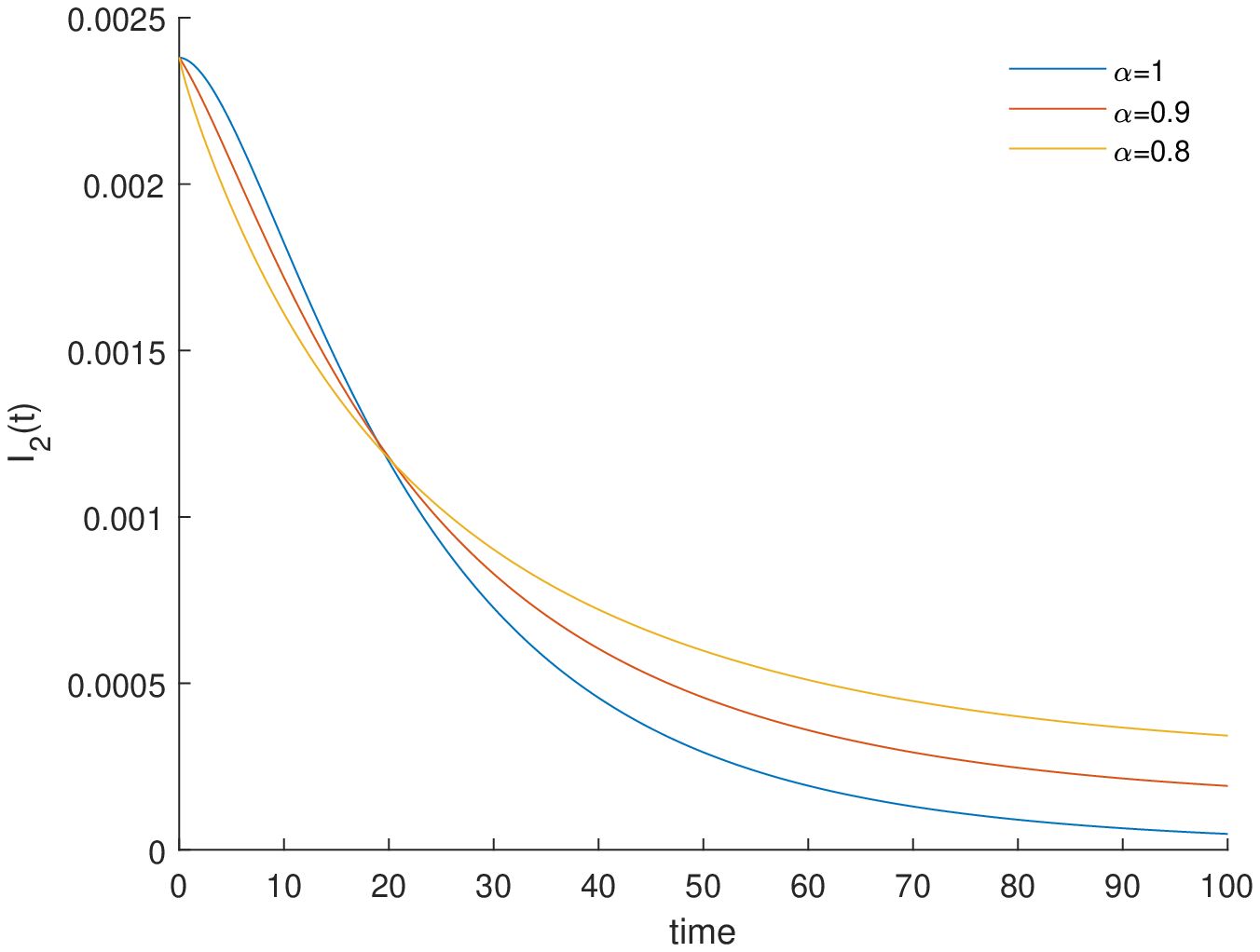}}}
\caption{Variables $I_1$ and $I_2$ of the FOCP \eqref{cost-functional}--\eqref{eq:ext:cont}
with values from Table~\ref{tab:param} and fractional order derivatives
$\alpha=1.0,$ $0.9$ and $0.8$. (\textbf{a}) Evolution of infected individuals of 1st community; 
(\textbf{b}) Evolution of infected individuals of 2nd community.}
\label{fig:infected_var:alphas}
\end{figure}

\begin{figure}[H]
\subfloat[]{%
\resizebox*{6.1cm}{!}{\includegraphics{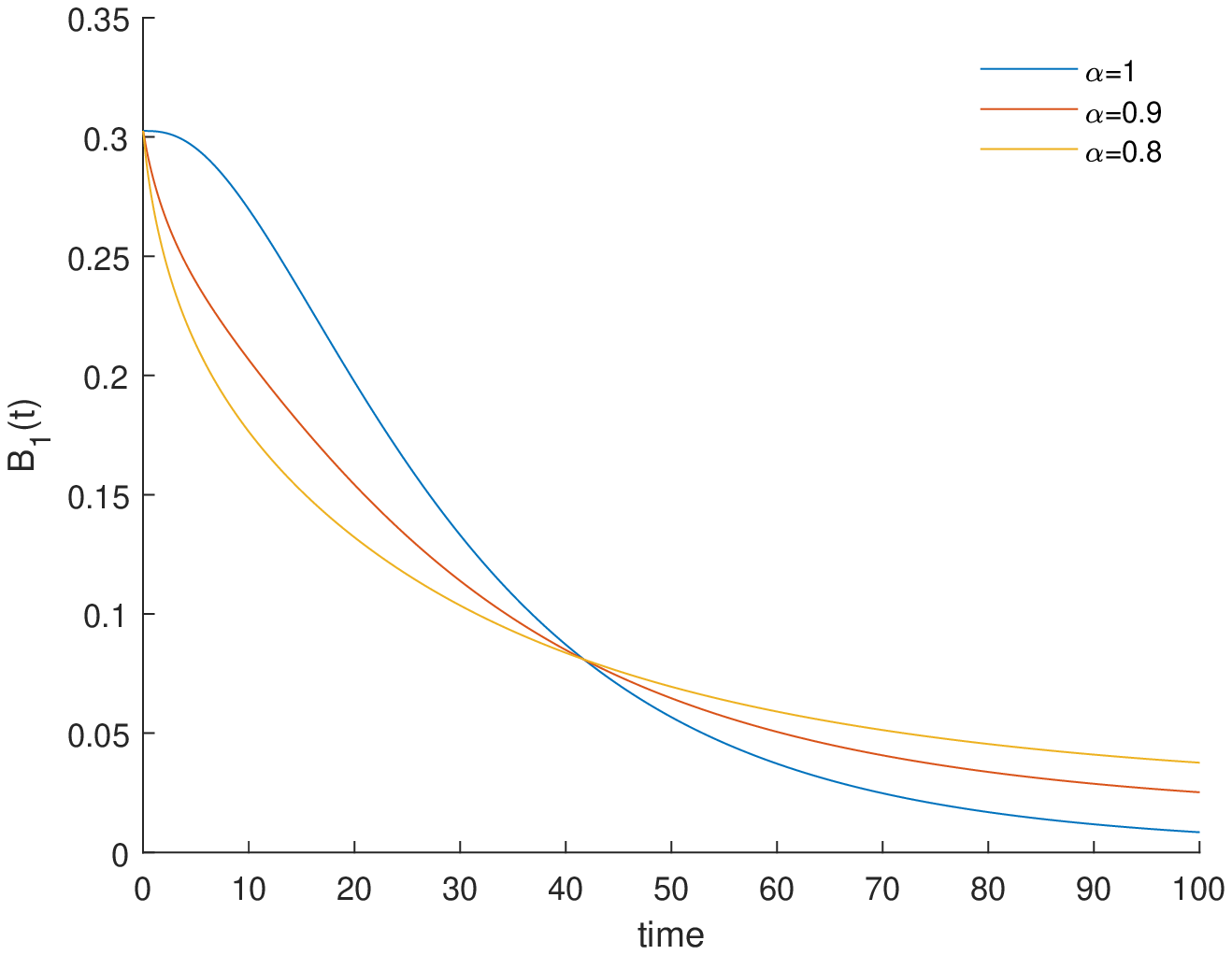}}}\hspace{4pt}
\subfloat[]{%
\resizebox*{6.1cm}{!}{\includegraphics{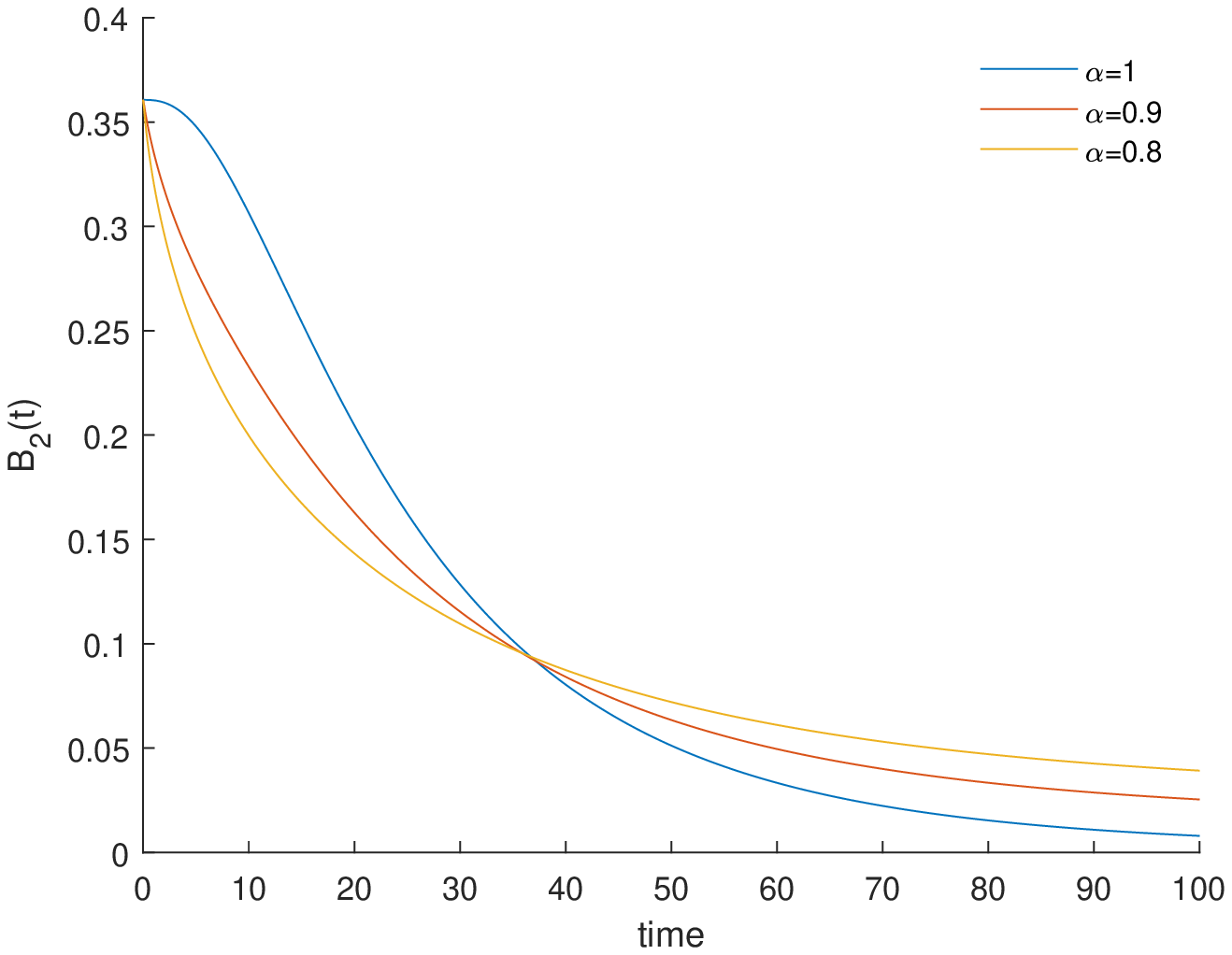}}}
\caption{Variables $B_1$ and $B_2$ of the FOCP
\eqref{cost-functional}--\eqref{eq:ext:cont}
with values from Table~\ref{tab:param} and fractional order derivatives
$\alpha=1.0$, $0.9$ and $0.8$. (\textbf{a}) Variation of \emph{vibrio} 
population in 1st community; (\textbf{b}) Variation of \emph{vibrio} population in 2nd community.}
\label{fig:vibrio_var:alphas}
\end{figure}

\begin{figure}[H]
\subfloat[]{%
\resizebox*{6.1cm}{!}{\includegraphics{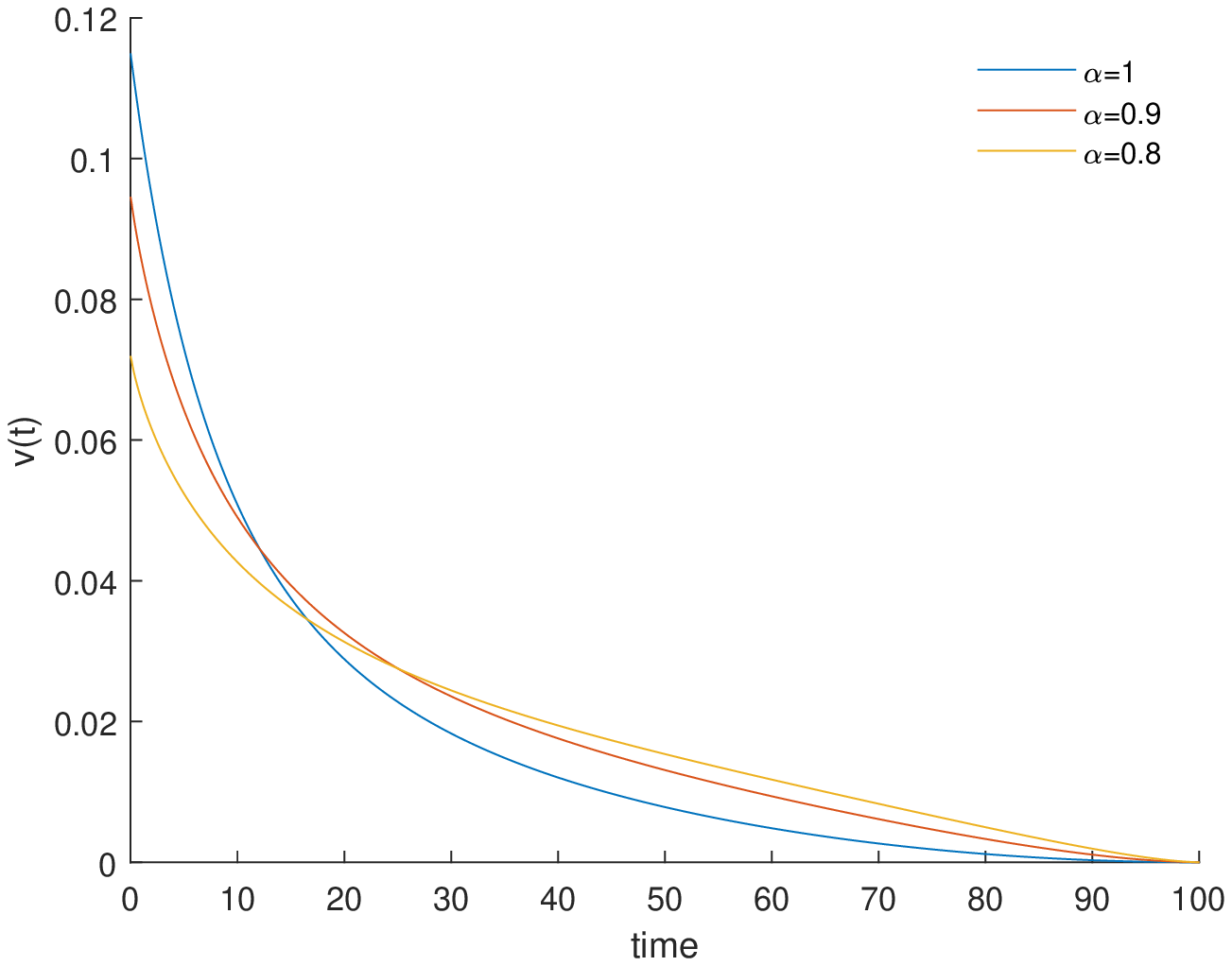}}}\hspace{4pt}
\subfloat[]{%
\resizebox*{6.1cm}{!}{\includegraphics{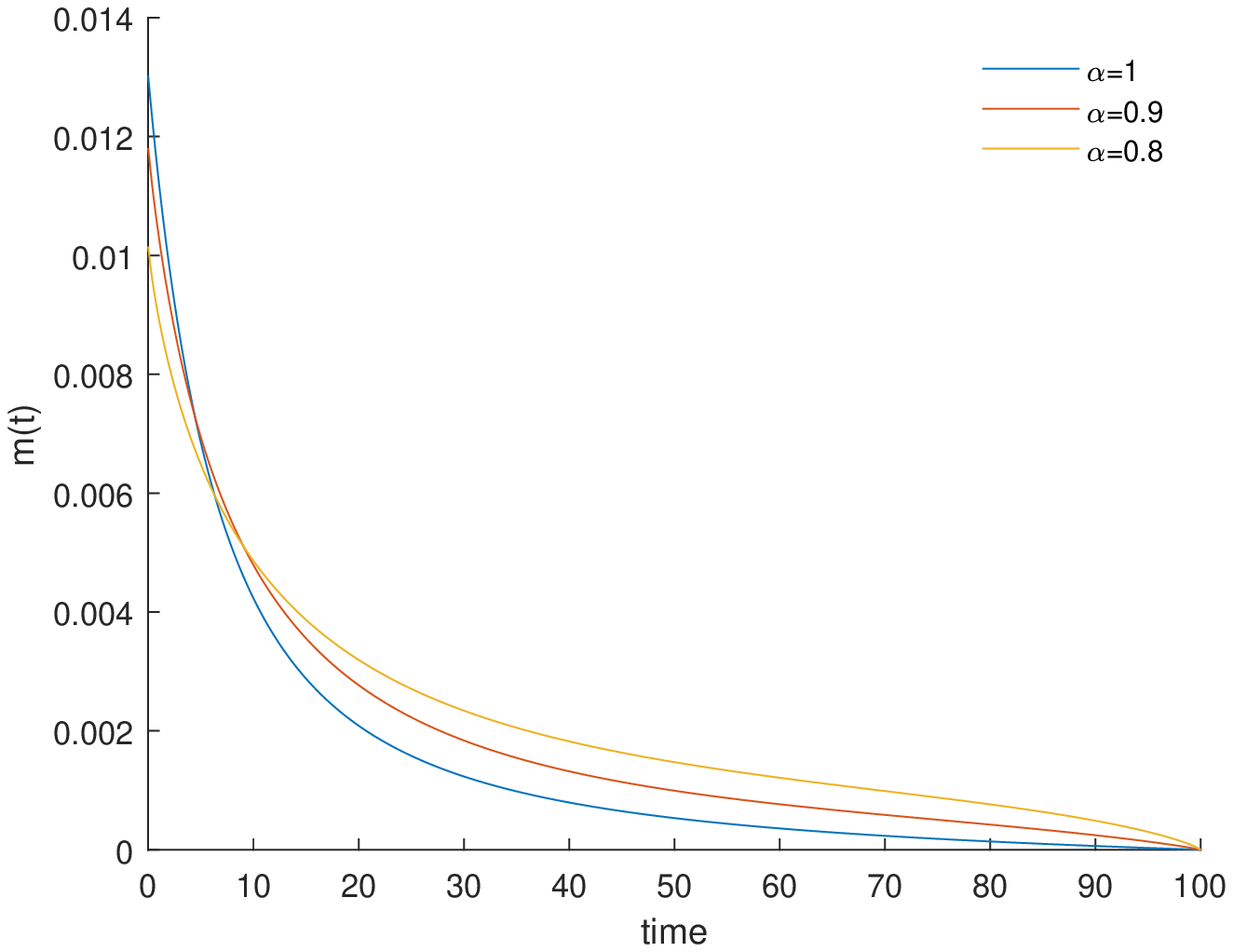}}}
\caption{Control variables $v$ and $m$ of the FOCP
\eqref{cost-functional}--\eqref{eq:ext:cont}
with values from Table~\ref{tab:param} and fractional order derivatives
$\alpha=1.0$, $0.9$ and $0.8$. (\textbf{a}) Contour of vaccination, $v$; 
(\textbf{b}) Evolution of control of hygiene, $m$.}
\label{fig:v_m_var:alphas}
\end{figure}

Without loss of generality, we consider the fractional order derivatives
$\alpha=1.0$, $0.9$ and $0.8$. In~Figures~\ref{fig:infected_var:alphas}--\ref{fig:v_m_var:alphas}, 
we find the solutions of the fractional optimal control problem for those values of $\alpha$.
We can see that a change in the value of $\alpha$ corresponds to significant
variations of the state and control variables. Beyond~those values of
$\alpha$, others values were also tested, but~the
results do not changed~qualitatively.

\begin{figure}[H]
\subfloat{%
\resizebox*{7cm}{!}{\includegraphics{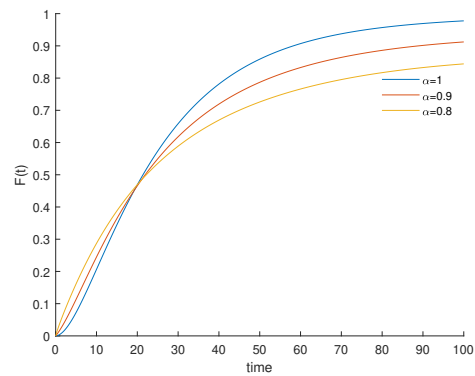}}}
\caption{Evolution of the efficacy function \eqref{efficacy_function}
for the FOCP \eqref{cost-functional}--\eqref{eq:ext:cont} with values from
Table~\ref{tab:param} and fractional order derivatives
$\alpha=1.0,$ $0.9$ and $0.8$.}
\label{fig:F_var:alphas}
\end{figure}

The efficacy function~\cite{rodrigues2014cost},
exhibited in Figure~\ref{fig:F_var:alphas}, is defined as
\begin{equation}
\label{efficacy_function}
F(t)=\frac{i(0)-i^*(t)}{i(0)}
=1-\frac{i^*(t)}{i(0)},
\end{equation}
where $i^*(t)=I_1^*(t)+I_2^*(t)$ is the optimal solution associated
with the fractional optimal control problem and $i(0)=I_1(0)+I_2(0)$ is the
correspondent initial condition. This function measures the proportional
variation  in the number of infected individuals, of~both communities,
after the application of the control measures, $\{v^*,m^*\}$, by~comparing
the number of infectious individuals at time $t$ with its initial value $i(0)$.
We observe that the graphic of $F(t)$ exhibits the inverse tendency of
infected individuals curves, growing and reaching the maximum
at the end of the time~interval.

To assess the cost and the effectiveness of the proposed fractional
control measure during the intervention period,
some summary measures are presented.
The total cases averted by the intervention during
the time period $t_f$ is defined
in~\cite{rodrigues2014cost} by
\begin{equation}
\label{eq:A}
AV=t_f i(0)-\int_0^{t_f}i^*(t)~dt,
\end{equation}
where $i^*(t)$ is the optimal solution associated with
the fractional optimal controls and $i(0)$ is the correspondent
initial condition. Note that this initial condition is obtained
as the equilibrium proportion $\overline{i}$ of systems
\eqref{Cholera_model:I}--\eqref{Cholera_model:II},
which is independent of time, so that
$t_f i(0)=\displaystyle \int_0^{t_f}\overline{i}~dt$
represents the total infectious cases over a given period of $t_f$ days.

Effectiveness is defined as the proportion of cases averted
on the total possible cases under no intervention
~\cite{rodrigues2014cost}:
\begin{equation}
\label{eq:F}
\overline{F}=\frac{AV}{i(0) t_f}
=1-\frac{\displaystyle \int_0^{t_f}i^*(t)~dt}{i(0) t_f}.
\end{equation}

The total cost associated with the intervention
is defined in~\cite{rodrigues2014cost} by
\begin{equation}
\label{eq:TCI}
TC=\int_0^{t_f}( C_1\,v^*(t)s^*(t)+C_2\,m^*(t)i^*(t))dt,
\end{equation}
where $s^*(t)=S_1^*(t)+S_2^*(t)$ and $C_i$ corresponds to the
per person unit cost of the two possible interventions:
(i) vaccination at any time $t$ of susceptible individuals ($C_1$);
and (ii) infected individuals practising proper 
hygiene ($C_2$).~Following~\cite{okosun2013optimal,rodrigues2014cost},
the average cost-effectiveness ratio is given by
\begin{equation}
\label{eq:ACER}
ACER=\frac{TC}{AV}.
\end{equation}

The cost-effectiveness measures are summarised in Table~\ref{tab:efficacy}.~The 
results show the effectiveness of controls to reduce cholera infectious
individuals and the leadership in doing so by the classical model ($\alpha=1$).
\begin{specialtable}[H]
\caption{Summary of cost-effectiveness measures for classical and fractional
($0<\alpha <1$) cholera disease optimal control problems. Parameters according
to Tables~\ref{tab:param} and~\ref{tab:solinit} with $C_1=C_2=1.$}\label{tab:efficacy}
\begin{tabular*}{\hsize}{@{}@{\extracolsep{\fill}}ccccc@{}}
\toprule
\boldmath{$\alpha$} & \boldmath{$AV$}  & \boldmath{$TC$} 
& \boldmath{$ACER$} & \boldmath{$\overline{F}$}  \\[1mm] \midrule
1.0  &  0.316716  &  0.90049   &  2.84322  &  0.723582 \\
0.9  &  0.297175  &  1.17595  &   3.95708  &  0.678938 \\
0.8  &  0.280311  &  1.35978  &   4.85099  &  0.640408 \\ \bottomrule
\end{tabular*}
\end{specialtable}


\unskip


\unskip


\unskip


\unskip


\subsubsection{Optimal Control Strategies and Cost-Effectiveness~Analysis}

Now, we analyse the cost-effectiveness of alternative combinations of two
possible  control~measures:
\begin{itemize}
\item strategy $A$---implementing the vaccination control, $v$;
\item strategy $B$---implementing proper hygiene control, $m$;
\item strategy $C$---implementing both controls, $v$ and $m$.
\end{itemize}

To analyse the cost-effectiveness of the three alternative strategies
$A$, $B$, and $C$, we use the Incremental Cost-Effectiveness Ratio (ICER)
\cite{okosun2013optimal}.~This ratio is used to compare the differences
between costs and health outcomes of two alternative intervention
strategies that compete for the same resources, being often described
as the additional cost per additional health outcome.~We start by ranking 
the strategies in order of increasing effectiveness,
assessed by the total averted cases $AV$, defined in \eqref{eq:A}.

The ICER for the classical model ($\alpha=1$), is calculated as follows:
\begin{equation*}
\begin{split}
\mbox{ICER($B$)}=& \mbox{ACER}\mbox{($B$)}=0.216276,\\[1mm]
\mbox{ICER($A$)}=& \frac{0.900865-0.0084106}{0.316411-0.038888}=3.2157877,\\[1mm]
\mbox{ICER($C$)}=& \frac{1.900494-0.900865}{0.316716-0.316411}=-1.216393.
\end{split}
\end{equation*}

Results are shown in Table~\ref{tab:effectiveness}.~Strategy $A$ is the most costly one, 
so we exclude this strategy from the set of alternatives.~We align the remaining strategies 
by increasing effectiveness ($AV$) and recalculate the ICER: 
$ICER(B)= ACER(B)=0.216276$ and $ICER(C)=3.210922$. Hence, we conclude 
that strategy $B$ (implementing only the control of hygiene measure, $m$)
is the most cost-effective~strategy.
\begin{specialtable}[H]
\caption{Incremental cost-effectiveness ratio for strategies $A$, $B$, and $C$.
Parameters according to Tables~\ref{tab:param} and~\ref{tab:solinit}
with $C_1=C_2=1$ and $\alpha=1$.}
\label{tab:effectiveness}
\begin{tabular*}{\hsize}{@{}@{\extracolsep{\fill}}ccccc@{}}
\toprule
\textbf{Strategies} & \boldmath{$AV$}  & \boldmath{$TC$} 
& \boldmath{$ACER$} & \textbf{ICER}  \\[1mm] \midrule
$B$  & 0.038888 & 0.0084106 & 0.216276 &  0.216276\\
$A$  & 0.316411 & 0.900865  & 2.84713  & 3.2157877\\
$C$  & 0.316716 & 0.900494  & 2.84322  & $-1.216393$\\
\bottomrule
\end{tabular*}
\end{specialtable}
The ICER was also computed for the fractional model, considering the same strategies,
with the three derivative orders used previously: $\alpha=0.9$, $\alpha=0.8$ and $\alpha=0.68$.
Even in those cases, we obtain the same conclusion: strategy $B$ is the most~cost-effective.

Comparing the cost-effectiveness of strategy $B$ for above derivative orders,
we start by getting the values of ICER presented in Table~\ref{tab:effectiveness2}.
Once the values of Total Cost (TC) diminish with the decrease in the derivative order,
proceeding as above, when the values of ICER were computed for the classical $\alpha=1$ model,
we eliminate successively the scheme with the highest TC. Consequently, the~``strategies''
$\alpha=1.0$, $\alpha=0.9$, and $\alpha=0.8$ are excluded, by~this~order.

Our conclusion is: the most cost-effective scheme is strategy $B$ with $\alpha=0.68$.

\begin{specialtable}[H]
\caption{Incremental cost-effectiveness ratio
for strategy $B$ and several derivative orders.
Same conditions as Table~\ref{tab:effectiveness}.}
\label{tab:effectiveness2}
\begin{tabular*}{\hsize}{@{}@{\extracolsep{\fill}}cccc@{}}
\toprule
\boldmath{$\alpha$} & \boldmath{$AV$}  
& \boldmath{$TC$} & \textbf{ICER}  \\[1mm] \midrule
1.0    & 0.038888 & 0.0084106 &  0.216276\\
0.9    & 0.147496 &  0.003549& $-22.3419$ \\
0.8    & 0.203782 &  0.001705& $-30.5191$  \\
0.68   & 0.237211 & 0.000845 & $-38.859$ \\
\bottomrule
\end{tabular*}
\end{specialtable}
\unskip


\subsubsection{The Variable-Order FractInt~System}
\label{subsec:fractint}

The solution of the Fractional Optimal Control Problem, exhibited
in Section~\ref{subsec:numres} for some derivative order values,
evidences that the fractional-order model can be more  effective
in part of the time interval---as shown in Figure~\ref{fig:F_var:alphas}---but 
classical model is more effective if we consider all the interval,
as shown in Table~\ref{tab:efficacy}. In~this section, we consider that
the derivative order varies along the interval, being fractional or classical
when more advantageous. According with Figure~\ref{fig:F_var:alphas},
$\alpha$ should be fractional at beginning and become classical/one after a certain time.
Such class of Variable-Order Fractional (VOF) systems~\cite{almeida2019variable}
are baptised here as \emph{FractInt}. The~derivative order of the
\emph{FractInt} system varies according with
\begin{equation}
\label{alpha_t}
\alpha(t)=
\begin{cases}
\alpha_0 &\mbox{ if } 0\leqslant t\leqslant t',\\
1 &\mbox{ if } t'<t\leqslant 100,
\end{cases}
\end{equation}
where $0<\alpha_0<1$. In~practice, we noticed that the resulting system
is as more effective as smaller the value of $\alpha_0$. In~view of
an endemic scenario, we consider in our simulations $\alpha_0=0.68$,
which the lowest value that can guarantee it.
With respect to the switching time, the~value considered is $t'=7$, which is the one
to which it corresponds the maximum value of~efficiency.

The \emph{FractInt} system is numerically solved with the procedure described above,
at the beginning of Section~\ref{subsec:numres}. In~each iteration of the procedure,
the predict-evaluate-correct-evaluate method is applied successively to each one
of the two initial value problems---associated to the state system and to the adjoint system
of the FOCP \eqref{cost-functional}--\eqref{eq:ext:cont}---to perform the integration
of them in two steps. These steps correspond to the two branches of the derivative order function,
defined by \eqref{alpha_t}. In~such process, the~final solution of first step
(first branch) is the initial solution of the second~step.

\begin{figure}[H]
\subfloat[]{%
\resizebox*{6.1cm}{!}{\includegraphics{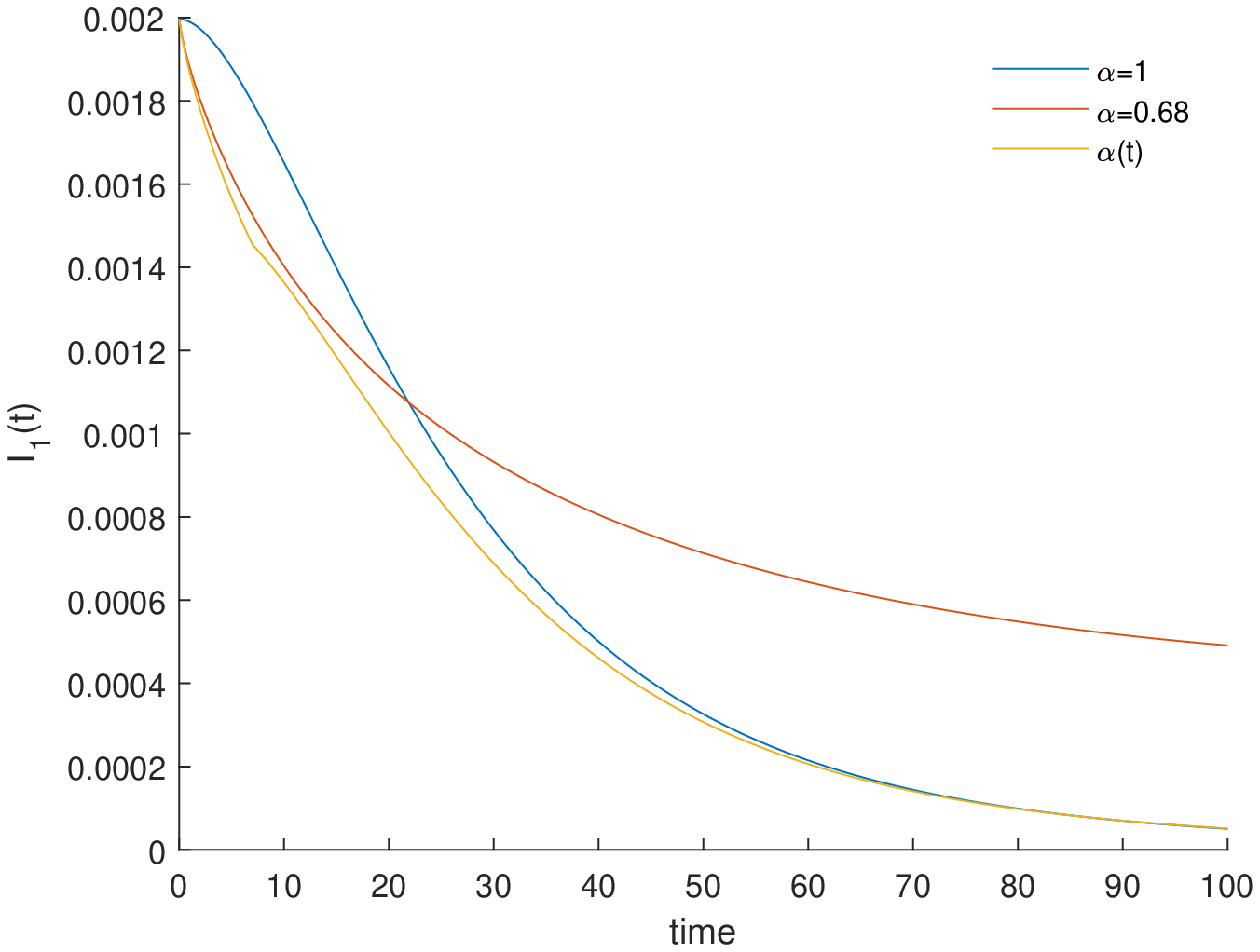}}}\hspace{4pt}
\subfloat[]{%
\resizebox*{6.1cm}{!}{\includegraphics{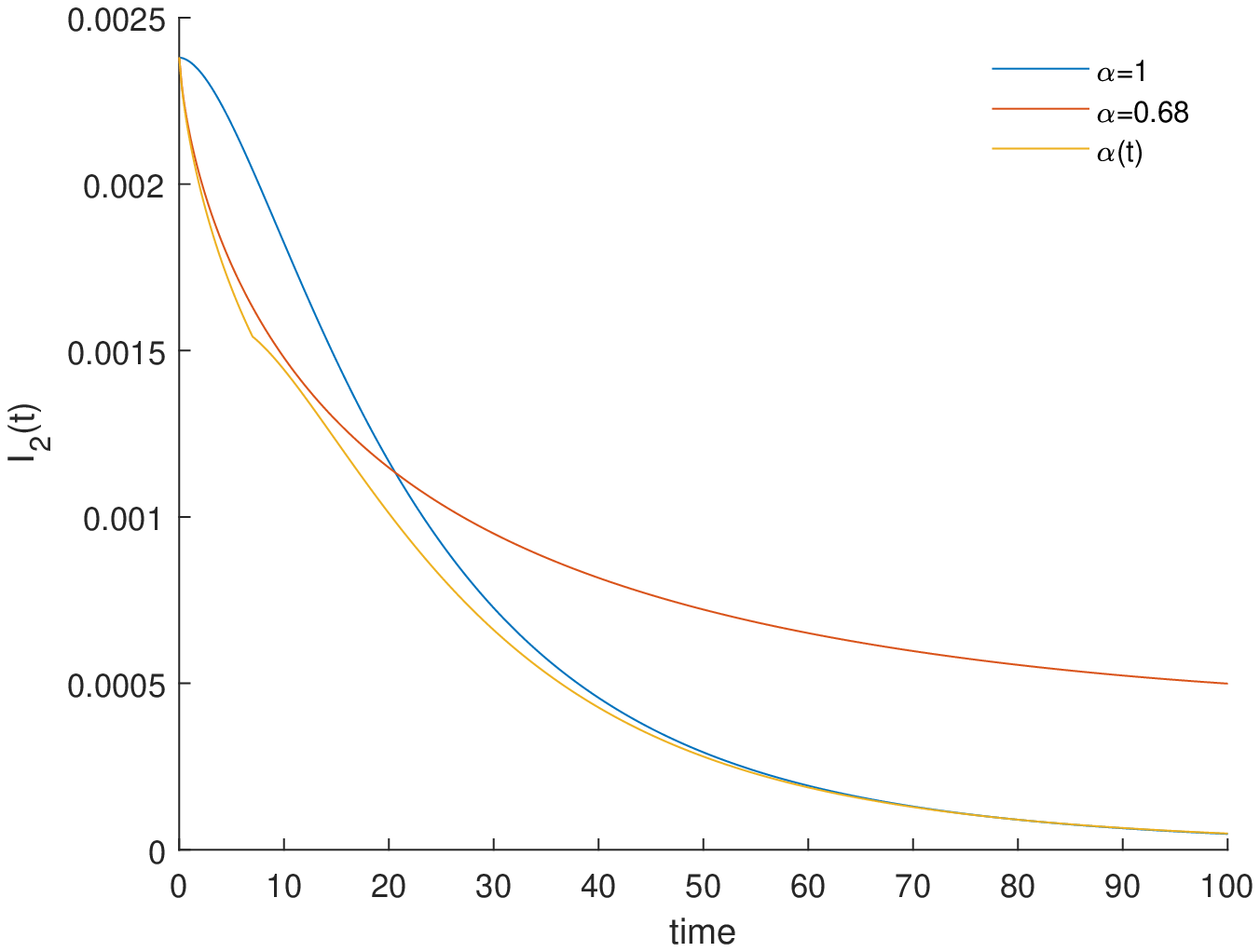}}}
\caption{Comparison of variables $I_1$ and $I_2$ of the
\emph{FractInt} system with the ones of the FOCP
\eqref{cost-functional}--\eqref{eq:ext:cont} with $\alpha=1$ and
$\alpha=0.68$, considering the parameter
values from Table~\ref{tab:param}. (\textbf{a}) Infected individuals 
of 1st community; (\textbf{b}) Infected individuals of 2nd community.}
\label{fig:infected_var:fractint}
\end{figure}

\begin{figure}[H]
\subfloat[]{%
\resizebox*{6.1cm}{!}{\includegraphics{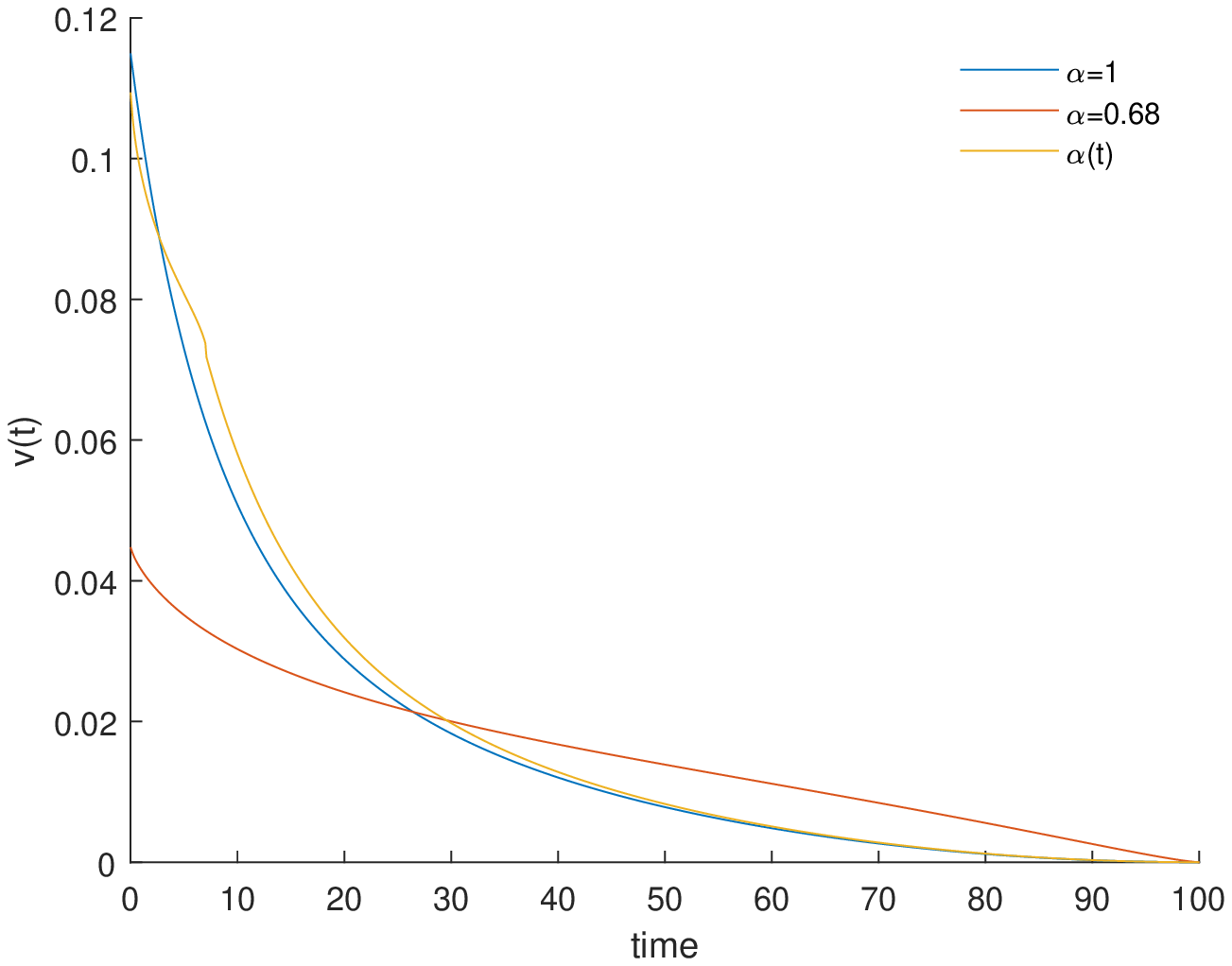}}}\hspace{4pt}
\subfloat[]{%
\resizebox*{6.1cm}{!}{\includegraphics{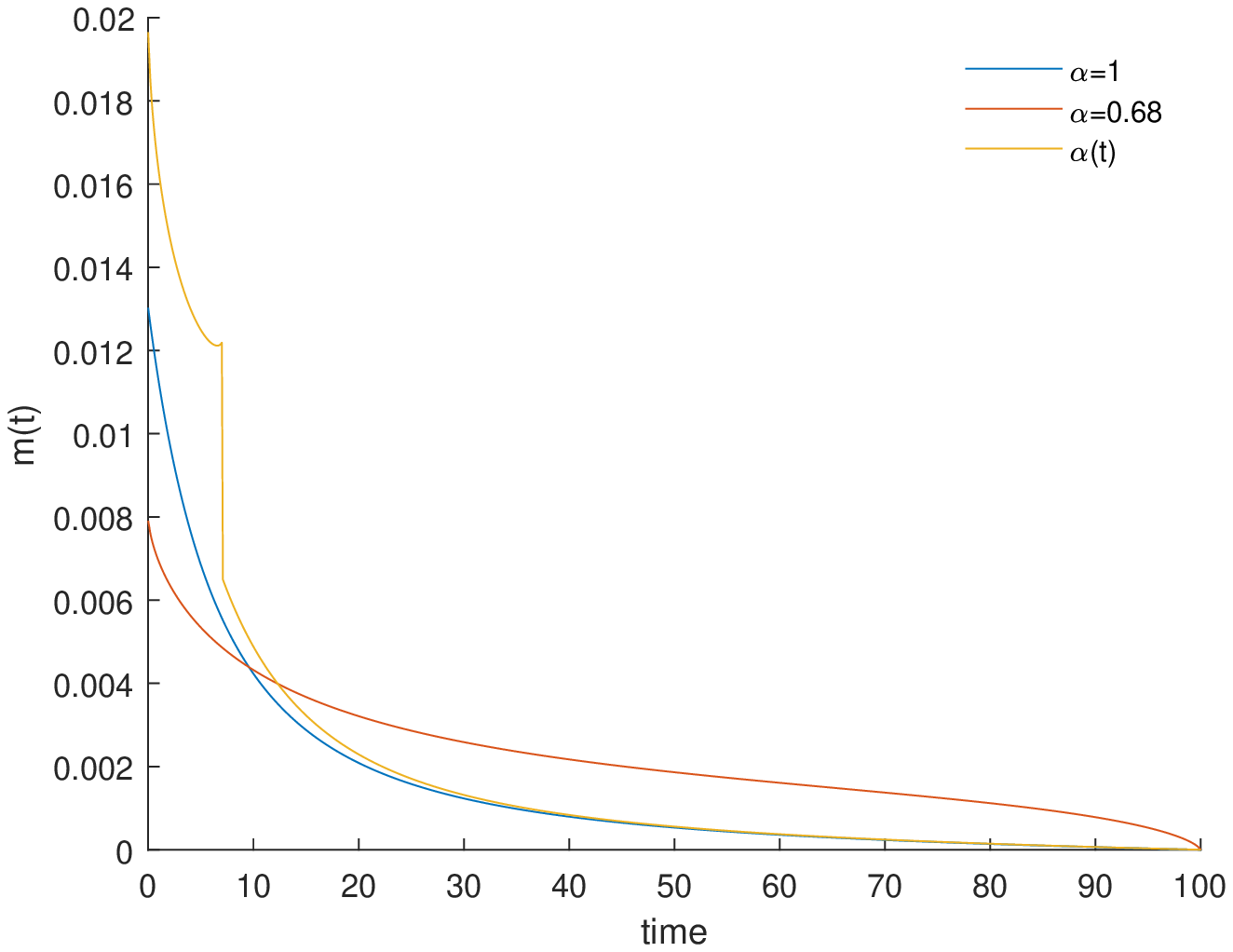}}}
\caption{Comparison of variables $v$ and $m$ of the \emph{FractInt} system with the
ones of the FOCP \eqref{cost-functional}--\eqref{eq:ext:cont} with $\alpha=1$
and $\alpha=0.68$, considering the parameter values from Table~\ref{tab:param}. 
(\textbf{a}) Contour of vaccination, $v$; (\textbf{b}) Evolution of the hygiene control, $m$.}
\label{fig:v_m_var:fractint}
\end{figure}

\begin{figure}[H]
\subfloat{%
\resizebox*{7cm}{!}{\includegraphics{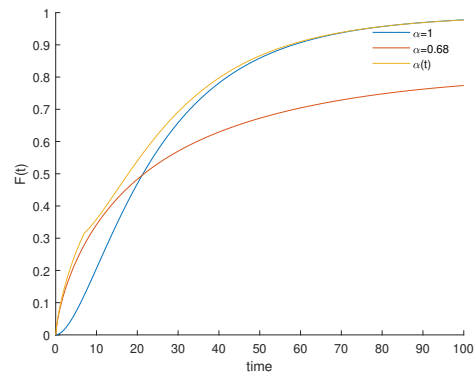}}}
\caption{Evolution of the efficacy function for the \emph{FractInt}
system and for the FOCP with $\alpha=1$ and $\alpha=0.68$,
considering the parameter values from Table~\ref{tab:param}.}
\label{fig:F_fractint:alphas}
\end{figure}

Solutions of the \emph{FractInt} system are reproduced in
Figures~\ref{fig:infected_var:fractint} and \ref{fig:v_m_var:fractint}
along with solutions of the classical and fractional models ($\alpha=1$
and $\alpha=0.68$, respectively).  We can see that the solutions of the
\emph{FracInt} system start to follow the solutions of the fractional model
and end by following the solutions of the classical model. This behaviour
can be observed in Figure~\ref{fig:F_fractint:alphas}, where
the efficacy function for these three models is~displayed.

The cost-effectiveness measures for the \emph{FractInt} system
are summarised in Table~\ref{tab:efficacy_fint}.

\begin{specialtable}[H]
\caption{Cost-effectiveness measures for the \emph{FractInt} system.
Parameters according to Tables~\ref{tab:param}
and~\ref{tab:solinit} with $C_1=C_2=1$.}
\label{tab:efficacy_fint}
\begin{tabular*}{\hsize}{@{}@{\extracolsep{\fill}}cccc@{}}
\toprule
\boldmath{$AV$}  & \boldmath{$TC$} & \boldmath{$ACER$} 
& \boldmath{$\overline{F}$}  \\[1mm] \midrule
0.332078 & 1.10967 &  3.3416 & 0.758679\\ \bottomrule
\end{tabular*}
\end{specialtable}

Our results show that it is effective to control cholera infection
through optimal control and that the \emph{FracInt} model
is more efficient than the classical one (cf. Table~\ref{tab:efficacy}
and \mbox{Figure~\ref{fig:F_fractint:alphas}}), being the most effective~model.


\unskip


\unskip


\unskip


\section{Conclusions}
\label{sec:conclusion}

Nowadays, Cholera infection is still an healthcare problem 
in many parts of the world, namely in impoverished regions 
where it has devastating effects. Its control is the goal 
of many studies in last years. In~this work,
a fractional-order mathematical model for
Cholera with two connected communities, proposed 
in~\cite{njagarah2018spatial}, is studied  and generalised.
The two communities in appreciation are distinct and therefore 
they could not be molten in one unique~community.

A sensitivity analysis of the model is done to show
the importance of estimation of parameters. A~fractional
optimal control (FOC) problem is then formulated and numerically~solved.

The relevance of studied controls is assessed using a cost-effectiveness 
analysis. Such analysis allow us to neglect the control $u$
(domestic water treatment) since it proved to be useless when 
used in combination with remaining~controls.

The numerical results show that the FOC system is more effective 
only in part of the time interval. Therefore, we propose a system 
where the derivative order varies along the time interval, being 
fractional or integer when more advantageous
in the control of infection.
Such variable-order fractional model, baptised here as \emph{FractInt},
shows to be the most effective in the control of the~disease.


\vspace{6pt}

\authorcontributions{Conceptualisation, S.R. and D.F.M.T.; 
methodology, S.R. and D.F.M.T.; software, S.R.; 
validation, S.R. and D.F.M.T.; formal analysis, S.R. and D.F.M.T.; 
investigation, S.R. and D.F.M.T.; 
writing---original draft preparation, S.R. and D.F.M.T.; 
writing---review and editing, S.R. and D.F.M.T.; 
visualisation, S.R. and D.F.M.T. All authors have read and agreed 
to the published version of the~manuscript.}

\funding{This research was funded by 
Funda\c c\~ao para a Ci\^encia e a Tecnologia 
(FCT, the~Portuguese Foundation for Science and Technology)	
through IT, Grant Number UIDB/50008/2020 (S.R.),
and CIDMA, Grant Number UIDB/04106/2020 (D.F.M.T.).}

\institutionalreview{Not~applicable.} 

\informedconsent{Not applicable.} 

\dataavailability{Not applicable.} 

\acknowledgments{The authors would like to thank John B. H. Njagarah,
from Botswana International University of Science and Technology,
for having provided most of the values of parameters
used in this work, which are presented in Table~\ref{tab:param}.}

\conflictsofinterest{The authors declare no conflicts of interest. 
The funders had no role in the design of the study; in the collection, 
analyses, or~interpretation of the data; in the writing of the manuscript; 
nor in the decision to publish the~results.} 


\end{paracol}

\reftitle{References}


\end{document}